\documentclass[11pt]{article}
\usepackage{amsfonts,latexsym}
\usepackage{amsmath}
\usepackage{amsthm}
\usepackage{amssymb}
\usepackage{mathrsfs}

\newtheorem{theorem}{\indent Theorem}[section]
\newtheorem{proposition}[theorem]{\indent Proposition}
\newtheorem{definition}[theorem]{\indent Definition}
\newtheorem{lemma}[theorem]{\indent Lemma}
\newtheorem{remark}[theorem]{\indent Remark}

\begin{document}

\title{Uniform attractors for non-autonomous wave equations\\ with nonlinear damping }
\author{Chunyou Sun, Daomin Cao\\
{\small Institute of Applied Mathematics, Chinese Academy of Sciences}\\
{\small Beijing, 100080, China} \\
{\small E-mails: cysun@amss.ac.cn;  dmcao@amt.ac.cn}\\
and \\
Jinqiao Duan\\
{\small Department of Applied Mathematics, Illinois Institute
of Technology} \\
{\small Chicago, IL 60616, USA} \\
{\small E-mail:  duan@iit.edu}}
\date{}

\renewcommand{\theequation}{\arabic{section}.\arabic{equation}}
\numberwithin{equation}{section}

\maketitle

\begin{abstract}
We consider dynamical behavior of   non-autonomous wave-type
evolutionary equations with nonlinear damping, critical
nonlinearity, and time-dependent  external forcing which is
translation bounded but not translation compact (i.e., external
forcing is not necessarily time-periodic, quasi-periodic or almost
periodic). A sufficient and necessary condition for the existence
of uniform attractors is established using the concept of uniform
asymptotic compactness. The required compactness for the existence
of uniform attractors is then fulfilled by some new a priori
estimates for concrete wave type equations arising from
applications. The structure of   uniform attractors is obtained by
constructing a skew product flow on the extended phase space for
the norm-to-weak continuous process.

{\bf Keywords:} Non-autonomous systems; Wave equation; Nonlinear
damping; Critical exponent; Uniform attractor; Norm-to-weak
continuous process.
\end{abstract}


\footnote[0]{\hspace*{-7.4mm}
Date: June 26, 2006.\\
AMS Subject Classification: 35L05, 35B40, 35B41\\
This work was partly supported by the NSF Grants DMS-0209326 \&
DMS-0542450, and the Outstanding Overseas Chinese Scholars Fund of
the Chinese Academy of Sciences.}

\section{Introduction}
We consider the following non-autonomous wave equations with
nonlinear damping, on a bounded domain $\Omega$ in $\mathbb{R}^3$
with smooth boundary:
\begin{equation}\label{1.1}
u_{tt}+h(u_t)-\Delta u+f(u)=g(x,\,t)\quad x\in\Omega
\end{equation}
under the boundary condition
\begin{equation}\label{1.2}
u |_{\partial \Omega}=0,
\end{equation}
and initial conditions
\begin{equation}\label{1.3}
u(x,0)=u_ {0}(x),\quad u_ t(x,0)=v_ 0(x).
\end{equation}
Here $h$ is the nonlinear damping function, $ f$ is the
nonlinearity, and $g$ is a given external time-dependent forcing.

Equation \eqref{1.1} arises as an evolutionary mathematical model
in various systems, for example: (i) modeling a continuous
Josephson junction with specific $h, g$ and $ f $; (ii) modeling a
hybrid system of nonlinear waves and nerve conduct; and(iii)
modeling a phenomenon in quantum mechanics. A relevant physical
issue is to investigate the asymptotic dynamical behavior of these
mathematical models.

For the special, autonomous case of \eqref{1.1}, i.e., when $g$
does not depend on time $t$ explicitly, the solution operator
defines a flow or semigroup. The asymptotic behaviors of the
solutions have been studied extensively by using of the concept of
global attractors; see, for example, \cite{ACH,Ba,BV,CV,Ha,Te} for
the linear damping case, and \cite{CL1,CL2,CL3,Fe,GM,SYZ1} for the
nonlinear damping case.

In the general case of non-autonomous system \eqref{1.1}, the
solution operator does \emph{not} define  a flow or semigroup, but
a \emph{process}; see Section 2 and Section 5 below. A proper
extension of the notion of a global attractor for semigroups to
the case of processes is the so-called uniform attractor (see e.g.
\cite{Ha3,CV94,CV}). About the basic concepts of non-autonomous
dynamical systems, uniform attractors and processes, we refer to
\cite{Ha3,CV94,CV} for more details or see Section 2 below.

 The basic assumptions about nonlinear damping
$h$ and nonlinearity $f$ are as follows:
\begin{equation}\label{1.4}
h\in C^1(\mathbb{R}),\quad h(0)=0,\quad h~\text{\rm is strictly
increasing},
\end{equation}
\begin{equation}\label{1.5}
\liminf_{|s|\to \infty}h'(s)>0,
\end{equation}
\begin{equation}\label{1.6}
|h(s)|\leq C_1(1+|s|^p),
\end{equation}
where $p \in [1,\, 5)$; $f\in C^1(\mathbb{R})$ and satisfies
\begin{equation}\label{1.7}
|f'(s)|\leq C_2(1+|s|^{q}),
\end{equation}
\begin{equation}\label{1.8}
\liminf_{|s|\to \infty}\frac{f(s)}{s}>-\lambda_1,
\end{equation}
where $0\leqslant q\leqslant 2$ and $\lambda_1$ is the first
eigenvalue of $-\Delta$ in $H_0^1(\Omega)$, and these assumptions
are similar to that for autonomous system and come from
\cite{CL1,CL2,Fe} etc.

In this paper, we consider the non-autonomous system
\eqref{1.1}-\eqref{1.3} via the uniform attractors of the
corresponding family of processes $\{U_{\sigma}(t,\tau)\}$,
$\sigma\in \Sigma$, especially with: (i) the nonlinear damping
(i.e., $h$ is a nonlinear function), (ii) the nonlinearity $f(u)$
has critical exponent ($q=2$), and (iii) the external forcing
$g(x,t)$ is \emph{not}  translation  compact in
$L^2_{loc}(\mathbb{R};L^2(\Omega))$.

In Chepyzhov \& Vishik\cite{CV}, for the linear damping case
$h(v)=kv$ with  a constant $k>0$ and $q<2$ (subcritical), for
system \eqref{1.1}-\eqref{1.3}, the authors obtained the existence
of a bounded uniform absorbing set if $g$ is translation bounded,
and the existence of a uniformly attractor when $g$ is translation
compact (e.g., when $g$ is time-periodic, quasi-periodic or almost
periodic). Under the assumptions that $g$ and $\partial_t g$ are
both in the space of bounded continuous functions
$C_b(\mathbb{R},L^2(\Omega ))$ and $h$ satisfies the growth bounds
$0<\alpha \leqslant h'(s)\leqslant \beta<\infty$ for some
constants $\alpha$ and $\beta$, Zhou \& Wang\cite{ZW} have proved
the existence of kernel sections and obtained the estimation of
the Hausdorff dimension of the kernel sections.

For the {\it existence} of uniform attractors, as in autonomous
case, some kind of compactness of the family of processes is a key
ingredient. The corresponding compactness assumption in
\cite{Ha3,CV94,CV} is that the family of processes
$\{U_{\sigma}(t,\tau)\}$, $\sigma\in \Sigma$ has a compact uniformly
absorbing set. The number $q=2$ is called the critical exponent,
since the nonlinearity $f$ is not compact in this case, which is an
essential difficulty in studying the asymptotic behavior even for
the autonomous case \cite{ACH,Ba,BV,CL1,CL2,Fe,SYZ1} etc.

About the case of $1<p <5$ for the nonlinear damping exponent $p$,
as mentioned in Haraux\cite{Ha2}, even for the bounded dissipation,
it becomes much more difficult when $g$ depends on $t$, and the
characterizations of dynamics for this case are unknown to the
authors. Moreover, the nonlinearity of $h$ also brings some
difficulties for us to prove the compactness, for example, for the
autonomous linearly damped wave equations, Ball \cite{Ba} gives a
very nice method to verify the necessary asymptotic compactness,
so-called energy methods by many other authors, and then this method
is generalized and given a general abstract framework for its
applications by \cite{LS,MRW1,MRW2} and others to both autonomous
and non-autonomous cases. However, for our problem, due to the
nonlinear damping, it seems to be difficult to apply the method of
Ball \cite{Ba}.

 The purpose of this paper is to obtain the existence and structure of
 the compact, in the norm
topology of $H_0^1\times L^2$, uniform attractor when the external
forcing $g(x,t)$ is \emph{not} translation compact in
$L^2_{loc}(\mathbb{R};L^2(\Omega))$. For the existence of uniform
attractors, a main approach in \cite{CV94,CV} is by constructing
skew product flow on the extend phase space $X\times \Sigma$. They
require that the symbol space has some compactness so that the
skew product flow has some corresponding compactness, i.e., the
concept of translation compact functions (e.g., see
\cite{CV94,CV95a,CV}). Consequently, the compact uniform
attractors are obtained for the systems with symbols of compact
hulls, and the weakly compact uniform attractors for the systems
with symbols of weakly compact hulls. However, there are some
results show that one can obtain the compact uniform attractors
for the system with translation noncompact external forcing: by
generalizing the methods in \cite{MWZ}, the authors in \cite{LWZ}
obtain the existence of an uniform attractor for 2D Navier-Stokes
equation in bounded domain with a kind of translation noncompact
external forcing; in Zelik\cite{Ze1}, by use of a bootstrap
argument together with a sharp use of Gronwall-type lemmas, when
$h(v)\equiv kv$ and $g,\partial_tg\in
L^{\infty}(\mathbb{R};L^2(\Omega))$, the author obtains some
regularity estimates for the solutions of \eqref{1.1}, which
implies naturally the existence of an uniform attractor; see also
the results in Chepyzhov \& Vishik\cite{CV93}.

Furthermore, we consider the {\it structure} of the uniform
attractor by investigating the kernel sections of a process (see
\cite{CV94,CV} for more details).

Here, for system \eqref{1.1}-\eqref{1.3}, we further assume that
\begin{equation}\label{1.9}
g(\cdot,t)\in L^{\infty}(\mathbb{R};L^2(\Omega))
\end{equation}
and
\begin{equation}\label{1.10}
\partial_tg\in L_b^r(\mathbb{R};L^r(\Omega))~\text{with}~r>\frac
65,
\end{equation}
where the space $L_b^r(\mathbb{R};L^r(\Omega))$ of ``translation
bounded" functions will be defined in the beginning of the next
section. Roughly speaking, the two conditions \eqref{1.9}-
\eqref{1.10} mean that the external forcing $g$ is bounded in time
and its time-derivative $\partial_t g$ is translation bounded. It
is clear  that a function $g$ satisfies \eqref{1.9} and
\eqref{1.10} dose not need to be translation compact in
$L^2_{loc}(\mathbb{R};L^2(\Omega))$. Moreover, we remark that the
technical hypothesis \eqref{1.9}   is mainly for the existence of
a bounded uniformly absorbing set; see $Theorem$ \ref{t5.5} below
or \cite{Ha2} for more details.

It is interesting to note that if \eqref{1.10} is replaced by the
assumption that $g$ is translation compact (e.g., $g$ is a
periodic, quasi-periodic or almost periodic function in
$L^2_{loc}(\mathbb{R};L^2(\Omega))$), then our result on uniform
attractors (see Theorem \ref{t5.12}) in Section 5 still holds, but
the proof can be largely simplified; see
\eqref{5.28}-\eqref{5.29}, $Remark$ \ref{r5.13} and  $Remark$
\ref{r5.13a} below. At the same time, the method  in \cite{LWZ}
can not be applied to our problem as \eqref{1.1} is a hyperbolic
type of equation, and the decomposition  \& regularization method
in \cite{CV93,Ze1} appear also not applicable here   due to the
nonlinearity of $h$.

This paper is organized as follows. After introducing some basic
materials in Section 2,  we first present a criterion for the
existence of a compact uniform attractor in Section 3, using the
concept of uniform asymptotic compactness (different from the
corresponding concept in \cite{CV94,CV}) which is introduced by
Moise et al in \cite{MRW1} for the family of semi-processes. We
apply this concept to the family of processes; see $Definition$
\ref{d3.1} and $Theorem$ \ref{t3.4}. Then, we  investigate the
{\it structure} of the uniform attractor via kernel sections of a
process. In fact, we present     results on uniform attractors and
their decompositions into kernel sections for norm-to-weak
continuous processes (see $Definition$ \ref{d3.5}, $Theorem$
\ref{t3.8} and $Theorem$ \ref{t3.9}). Note that the norm-to-weak
continuity here is weaker than the usual {\it norm-to-norm} and
{\it weak-to-weak} continuities \cite{ZYS}.

In Section 4, partially inspired by the results in
\cite{CL1,CL2,CL3,Kh,SYZ2}, we present a simple method for
verifying the uniform asymptotic compactness for processes
generated by wave type evolutionary equations like \eqref{1.1};
see $Theorem$ \ref{t4.2}.

Finally, in Section 5, as applications to concrete wave type
evolutionary equations, we first prove the existence of compactly
uniform (w.r.t $\sigma\in \Sigma$) attractors when the external
forcing $g_0=\sigma_0$ satisfies \eqref{1.9} and \eqref{1.10}; see
$Theorem$ \ref{t5.12}; then by verifying the norm-to-weak
continuity we show that the uniform attractor w.r.t. initial time
$\tau$ of a process $\{U_{\sigma_0}(t,\tau)\}$ coincides with the
uniform attractor w.r.t. symbol $\sigma \in \Sigma'$ when the
external forcing $g_0=\sigma_0\in
W^{1,\,\infty}(\mathbb{R};L^2(\Omega))$; we further decompose this
uniform attractor into kernel sections; see $Lemma$ \ref{l5.18}
and $Theorem$ \ref{t5.19}.

\section{Preliminaries}

In this section, we recall some basic concepts about non-autonomous
systems. We refer to \cite{Ha3,CV94,CV} and the references therein
for more details.

Space of translation bounded functions in
$L^r_{loc}(\mathbb{R};L^k(\Omega))$, with $r, k \geqslant 1$:
\[
L^r_b(\mathbb{R};L^k(\Omega))=\{g\in
L^r_{loc}(\mathbb{R};L^k(\Omega)):~\sup_{t\in
\mathbb{R}}\int_t^{t+1}\left(\int_{\Omega}|g(x,t+s)|^kdx\right)^{\frac{r}{k}}ds<\infty\};
\]

Space of translation compact functions in
$L^2_{loc}(\mathbb{R};L^2(\Omega))$:
\begin{align*}
L^2_c(\mathbb{R};L^2(\Omega))= \Big\{&g\in
L^2_{loc}(\mathbb{R};L^2(\Omega)): \text{For any interval} \;
[t_1,t_2]\subset \mathbb{R},\\
& \{g(x,h+s):h\in \mathbb{R}\}| {[t_1,\,t_2]} \; \mbox{ is
precompact in} \; L^2(t_1,t_2;L^2(\Omega)) \Big\}.
\end{align*}

Let $X$ be a complete metric space and $\Sigma$ be a parameter set.

The operators $\{U_{\sigma}(t,\tau)\},\sigma\in \Sigma$ are said to
be a family of processes in $X$ with symbol space $\Sigma$ if for
any $\sigma\in \Sigma$
\begin{equation}\label{2.1}
U_{\sigma}(t,s)\circ U_{\sigma}(s,\tau)=U_{\sigma}(t,\tau),\quad
\forall ~t\geqslant s\geqslant \tau,~\tau\in \mathbb{R},
\end{equation}
\begin{equation}\label{2.2}
U_{\sigma}(\tau,\tau)=Id~(identity),\quad \forall
~\tau\in\mathbb{R}.
\end{equation}
Let $\{T(s)\}_{s\geqslant 0}$ be the translation semigroup on
$\Sigma$, we say that a family of processes
$\{U_{\sigma}(t,\tau)\}$, $\sigma\in \Sigma$ satisfies the
translation identity if
\begin{equation}\label{2.3}
U_{\sigma}(t+s,\tau+s)=U_{T(s)\sigma}(t,\tau),\quad \forall
~\sigma\in \Sigma,~t\geqslant \tau,~\tau\in \mathbb{R},~ s\geqslant
0,
\end{equation}
\begin{equation}\label{2.4}
T(s)\Sigma=\Sigma,\quad \forall ~ s\geqslant 0.
\end{equation}

By $\mathcal{B}(X)$ we denote the collection of the bounded sets of
$X$, and $\mathbb{R}^{\tau}=\{t\in \mathbb{R},\,t\geqslant \tau\}$.

\begin{definition}\cite{CV}
A bounded set $B_0\in \mathcal{B}(X)$ is said to be a bounded
uniformly (w.r.t. $\sigma\in \Sigma$) absorbing set for
$\{U_{\sigma}(t,\tau)\},\sigma\in \Sigma$ if for any
$\tau\in\mathbb{R}$ and $B \in \mathcal{B}(X)$ there exists $T_0 =
T_0(B,\tau)$ such that
$\bigcup_{\sigma\in\Sigma}U_{\sigma}(t;\tau)B\subset B_0$ for all
$t\geqslant T_0$.
\end{definition}

\begin{definition}\cite{CV}
A set $A\subset X$ is said to be uniformly (w.r.t $\sigma\in
\Sigma$) attracting for the family of processes
$\{U_{\sigma}(t,\tau)\},\sigma\in \Sigma$ if for any fixed $\tau\in
\mathbb{R}$ and any $B\in\mathcal{B}(X)$
\[
\lim_{t\to +\infty}\left(\sup_{\sigma\in \Sigma}dist (U_{\sigma}(t;
\tau)B;\, A )\right) = 0,
\]
here $dist(\cdot,\cdot)$ is the usual Hausdorff semidistance in $X$
between two sets.

In particular, a closed uniformly attracting set
$\mathscr{A}_{\Sigma}$ is said to be the uniform (w.r.t. $\sigma\in
\Sigma$) attractor of the family of processes
$\{U_{\sigma}(t,\tau)\},\sigma\in \Sigma$ if it is contained in any
closed uniformly attracting set (minimality property).
\end{definition}
Obviously, if the uniform (w.r.t. $\sigma\in \Sigma$) attractor
exists, it is unique.

In order to obtain the structure as well as the existence of the
uniform attractor, under the condition \eqref{2.3}-\eqref{2.4}, the
authors in \cite{CV} construct the skew product flow in $X\times
\Sigma$:
\begin{equation}\label{2.5}
S(t)(u,\sigma)=(U_{\sigma}(t,0)u,\,T(t)\sigma),\quad t\geqslant
0,~(u,\sigma)\in X\times \Sigma,
\end{equation}
and $\{S(t)\}_{t\geqslant 0}$ forms a semigroup on $X\times \Sigma$.

\section{Abstract results}

\subsection{Existence of the uniform attractor}

In this subsection, we present a criterion for the existence of
compact uniform  attractor using the concept of uniform  (w.r.t.
$\sigma\in\Sigma$) asymptotical  compactness, which is different
from the corresponding concept in \cite{CV94,CV}, and it is
introduced in Moise et al \cite{MRW1} for a family of
semi-processes. Now, we use  this concept to the family of
processes.

\begin{definition}\cite{MRW1}\label{d3.1}
A family of processes $\{U_{\sigma}(t,\tau)\}$, $\sigma\in \Sigma$
on a complete metric space $X$ is said to be uniformly (w.r.t.
$\sigma\in\Sigma$) asymptotically compact, if and only if for any
fixed $\tau\in \mathbb{R}$, bounded sequence
$\{u_n\}_{n=1}^{\infty}\subset X$,
$\{\sigma_n\}_{n=1}^{\infty}\subset \Sigma$ and any
$\{t_n\}_{n=1}^{\infty}\subset \mathbb{R}^{\tau}$ with $t_n\to
\infty$ as $n\to \infty$, sequence
$\{U_{\sigma_n}(t_n,\,\tau)u_n\}_{n=1}^{\infty}$ is precompact in
$X$.
\end{definition}

Similarly, define the uniform $\omega-$limit set of $B\subset X$ at
initial time $\tau$ by
\begin{equation}
\omega_{\tau,\Sigma}(B)=\bigcap_{t\geqslant\tau}
\overline{\bigcup_{\sigma\in\Sigma}\bigcup_{s\geqslant
t}U_\sigma(s,\tau)B},
\end{equation}
where $\overline{A}$ means the closure of $A$ in $X$.

Then, we have the following characterizations for the uniform
$\omega-$limit set, see \cite{CV,MRW1},
\begin{proposition}\label{p3.2}
For any bounded set $B\subset \mathcal{B}(X)$, $u\in
\omega_{\tau,\Sigma}(B)$ if and only if there exist
$\{u_n\}_{n=1}^{\infty}\subset B$,
$\{\sigma_n\}_{n=1}^{\infty}\subset \Sigma$ and
$\{t_n\}_{n=1}^{\infty}\subset \mathbb{R}^{\tau}$ with $t_n\to
\infty$ as $n\to \infty$ such that $U_{\sigma_n}(t_n,\,\tau)u_n\to
u$.
\end{proposition}

In the following, similar to \cite{LWZ}, we give some
characterizations for the uniform (w.r.t. $\sigma\in\Sigma$)
asymptotically compact processes.
\begin{lemma}\label{l3.3}
Let $\{U_\sigma(t,\tau )\} $, $\sigma \in \Sigma $ be a family of
uniform (w.r.t. $\sigma \in \Sigma$) asymptotically compact
processes on a complete metric space $X$, then for any
$\tau\in\mathbb{R}$ and any nonempty set $B\in\mathcal B(X)$, we
have
\begin{enumerate}
\item[(i)] $\omega_{\tau,\,\Sigma}(B)$ is nonempty and compact in
$X$;

\item[(ii)] $\lim\limits_{t\rightarrow +\infty }
\sup\limits_{\sigma \in \Sigma }\text{dist}\left( U_\sigma \left(
t,\tau \right)B,\omega_{\tau,\Sigma}(B)\right) =0$;

\item[(iii)]
if $Y$ is closed and uniform (w.r.t. $\sigma \in \Sigma$) attracts
$B$, then $\omega_{\tau,\Sigma}(B)\subset Y$;

\hspace{-0.8 cm} furthermore, if $\{ U_\sigma ( t,\tau )\} $,
$\sigma \in \Sigma $ satisfies the translation identity
\eqref{2.3}-\eqref{2.4}, then
\item[(iv)] $\omega_{\tau,\Sigma}(B)\equiv\omega_{0,\Sigma}(B)$, that is, $\omega_{\tau,\Sigma}(B)$
is independent of $\tau\in\mathbb{R}$.
\end{enumerate}
\end{lemma}
{\bf Proof.} (i) For any fixed $\tau \in \mathbb{R}$, then for any
$t_n\in \mathbb{R}^{\tau}$, $t_n\to \infty$, $\sigma_n\in\Sigma$ and
$x_n\in B$, by the definition of uniform (w.r.t. $\sigma\in\Sigma$)
asymptotic compactness we know that
$\{U_{\sigma_n}(t_n,\,\tau)x_n\}_{n=1}^{\infty}$ is precompact in
$X$, without loss of generality, we assume that
\[
U_{\sigma_n}(t_n,\,\tau)x_n \to y.
\]
Then by the definition of $\omega-$limit set we know that $y\in
\omega_{\tau,\,\Sigma}(B)$, which implies that
$\omega_{\tau,\,\Sigma}(B)$ is nonempty.

For any $y_m\in \omega_{\tau,\,\Sigma}(B)$, $m=1,2,\cdots$, we will
show that $\{y_m\}_{m=1}^{\infty}$ is precompact in $X$. By the
definition, for each $m\in \mathbb{N}$, there exist $t_m\in
\mathbb{R}^{\tau}$, $t_m \geqslant m$, $\sigma_m\in\Sigma$ and
$x_m\in B$ such that
\[
\rho(U_{\sigma_m}(t_m,\,\tau)x_m,\, y_m) \leqslant \frac{1}{m},
\]
where $\rho(\cdot,\cdot)$ is the metric on $X$.

Therefore, by the assumption of uniform (w.r.t. $\sigma\in\Sigma$)
asymptotic compactness again, we have
$\{U_{\sigma_m}(t_m,\,\tau)x_m\}_{m=1}^{\infty}$ is precompact in
$X$, and without loss of generality, we assume that
$\{U_{\sigma_m}(t_m,\,\tau)x_m\}_{m=1}^{\infty}$ is a Cauchy
sequence in $X$. Then, from
\begin{align*}
&\rho(y_n,\,y_m) \nonumber \\
&\leqslant
\rho(y_n,\,U_{\sigma_n}(t_n,\,\tau)x_n)+\rho(U_{\sigma_n}(t_n,\,\tau)x_n,\,U_{\sigma_m}(t_m,\,\tau)x_m)
+\rho(U_{\sigma_m}(t_m,\,\tau)x_m,\,y_m) \nonumber \\
&\leqslant \frac{1}{n} +
\rho(U_{\sigma_n}(t_n,\,\tau)x_n,\,U_{\sigma_m}(t_m,\,\tau)x_m)
+\frac{1}{m},
\end{align*}
we know that $\{y_m\}_{m=1}^{\infty}$ is also a Cauchy sequence in
$X$. Moreover, from the definition we obviously have that
$\omega_{\tau,\,\Sigma}(B)$ is closed in $X$.

Hence, $\omega_{\tau,\,\Sigma}(B)$ is compact in $X$.

(ii) If (ii) is not true, then there exist $\varepsilon_0>0$,
$\sigma_n\in\Sigma$, $x_n\in B$ and $t_n\in \mathbb{R}^{\tau}$ with
$t_n \geqslant n$, such that
\[
\mathrm{dist}(U_{\sigma_n}(t_n,\,\tau)x_n,\,
\omega_{\tau,\,\Sigma}(B)) \geqslant \varepsilon_0 \quad
~n=1,2,\cdots.
\]
However, the uniform (w.r.t. $\sigma\in\Sigma$) asymptotic
compactness implies that
$\{U_{\sigma_n}(t_n,\,\tau)x_n\}_{n=1}^{\infty}$ is precompact in
$X$, that is, $\{U_{\sigma_n}(t_n,\,\tau)x_n\}_{n=1}^{\infty}$ has a
convergent subsequence which converges to some point of
$\omega_{\tau,\,\Sigma}(B)$. This is a contradiction.

(iii) $\forall~y\in \omega_{\tau,\,\Sigma}(B)$. Then there are
$\sigma_n\in\Sigma$, $x_n\in B$ and $t_n\in \mathbb{R}^{\tau}$ with
$t_n \to \infty$ such that $U_{\sigma_n}(t_n,\,\tau)x_n \to y$. From
the assumption that $Y$ uniform attracts $B$, obviously, we have
\[
\mathrm{dist}(U_{\sigma_n}(t_n,\,\tau)x_n,\,Y)\to 0\quad
\text{as}~n\to \infty.
\]
At the same time, the closeness of $Y$ implies $y\in Y$. Hence,
$\omega_{\tau,\,\Sigma}(B) \subset Y$.

(iv) For any fixed $\tau \in \mathbb{R}$ and $\sigma\in \Sigma$,
from the translation identity \eqref{2.3} we know (e.g., see
\cite{LWZ,MRW1}) that for any $\tau_0\in \mathbb{R}$ there is a
$\sigma'\in \Sigma$ such that
\[
U_\sigma(t,\tau)=U_{\sigma'}(t-\tau+\tau_0,\tau_0),\quad \forall\,
t\geqslant \tau.
\]
Combining with \eqref{2.4}, we have that for any $t\geqslant \tau$,
\[
\bigcup_{\sigma\in\Sigma}\bigcup_{s\geqslant t}U_\sigma(s,\tau)B =
\bigcup_{\sigma\in\Sigma}\bigcup_{s\geqslant t}U_\sigma(s-\tau,0)B.
\]

Therefore, we have
\[
\omega_{\tau,\,\Sigma}(B)=\bigcap_{t\geqslant\tau}
\overline{\bigcup_{\sigma\in\Sigma}\bigcup_{s\geqslant
t}U_\sigma(s,\tau)B} = \bigcap_{t\geqslant 0}
\overline{\bigcup_{\sigma\in\Sigma}\bigcup_{s\geqslant
t}U_\sigma(s,0)B}=\omega_{0,\,\Sigma}(B).
\]
$\hfill \blacksquare$

\begin{theorem}\label{t3.4}
Let $X$ be a complete metric space, $\{U_\sigma( t,\tau )\}$,
$\sigma \in \Sigma $ be a family of processes on $X$ and satisfies
the translation identity \eqref{2.3}-\eqref{2.4}. Then, $\{U_\sigma
(t,\tau)\}$, $\sigma \in \Sigma $ has a compactly uniform (w.r.t.
$\sigma \in \Sigma$) attractor $\mathscr{A}_\Sigma$ in $X$ and
satisfies
\begin{equation*}
\mathscr{A}_\Sigma
=\omega_{0,\Sigma}\left(B_0\right)=\omega_{\tau,\Sigma}\left(B_0\right)
=\bigcup_{B\in\mathcal{B}(X)}\omega_{\tau,\Sigma}\left(B\right),
\quad\forall\, \tau\in \mathbb{R},
\end{equation*}
if and only if $\{U_\sigma (t,\tau )\}$, $\sigma \in \Sigma $
\begin{enumerate}
\item[i)] has a bounded uniformly (w.r.t. $\sigma \in
\Sigma$) absorbing set $B_0$;

\item[ii)] is uniformly (w.r.t. $\sigma \in \Sigma$)
asymptotically compact.
\end{enumerate}
\end{theorem}
{\bf Proof.} The necessity follows from the definition of uniform
(w.r.t $\sigma\in \Sigma$) attractor and the compactness of
$\mathscr{A}_{\Sigma}$.

Now we prove the sufficiency. For any fixed $\tau \in \mathbb{R}$
and any $B\in\mathcal{B}(X)$. We know that there is a $T=T(\tau,B)$
such that
\[
\bigcup_{\sigma\in \Sigma}\bigcup_{t\geqslant T}U_{\sigma}(t,\tau)B
\subset B_0.
\]
Combining with the equivalent characterization, $Proposition$
\ref{p3.2}, of $\omega-$limit set and $Lemma$ \ref{l3.3}, we have
\begin{equation}\label{3.2}
\omega_{\tau,\Sigma}\left(B\right) \subset
\omega_{\tau,\Sigma}\left(B_0\right)=\omega_{0,\Sigma}\left(B_0\right),
\end{equation}
and $\omega_{\tau,\Sigma}(B)$, of course,
$\omega_{0,\Sigma}\left(B_0\right)$ uniform (w.r.t. $\sigma \in
\Sigma$) attracts $B$.

Moreover, \eqref{3.2} implies that
$\bigcup_{B\in\mathcal{B}(X)}\omega_{\tau,\Sigma}\left(B\right)
\subset \omega_{0,\Sigma}\left(B_0\right)$, and from $B_0\in
\mathcal{B}(X)$ we obtain
$\bigcup_{B\in\mathcal{B}(X)}\omega_{\tau,\Sigma}\left(B\right) =
\omega_{0,\Sigma}\left(B_0\right)$.

The minimality and closeness follows immediately from (iii) of
$Lemma$ \ref{l3.3}, and the compactness follows from (i) of $Lemma$
\ref{l3.3}. $\hfill \blacksquare$

\subsection{Structure of the uniform attractor}

We describe the structure of the uniform attractor by means of its
kernel sections.

Hereafter, we assume that $X$ is a Banach space with norm
$\|\cdot\|_X$ and $\Sigma$ is a complete metric space with metric
$d(\cdot,\,\cdot)$.

Let $\{U(t,\tau)|t\geqslant \tau,~\tau\in
\mathbb{R}\}=\{U(t,\tau)\}$ be a process acting in a Banach space
$X$, and let $\mathcal{K}$ be the kernel of the process
$\{U(t,\tau)\}$. We recall (e.g., see \cite{CV}) that the kernel
$\mathcal{K}$ consists of all bounded complete trajectories of the
process, i.e.,
\[
\mathcal{K}=\{u(\cdot)~| \|u(t)\| X\leqslant
C_u,~U(t,\tau)u(\tau)=u(t),~\forall\,t\geqslant \tau,~\tau\in
\mathbb{R}\},
\]
and $\mathcal{K}(s)$ denotes the kernel section at a time moment
$s\in \mathbb{R}$:
\[
\mathcal{K}(s)=\{u(s)~|u(\cdot)\in \mathcal{K}\},\quad
\mathcal{K}(s)\subset X.
\]

As mentioned in \cite{CV}, since the invariance of the global
attractor of a semigroup is replaced by the minimality in the
definition of the uniform attractor of a family of processes, the
existence of uniform attractor dose not need any continuity for
the processes. However, in order to obtain the structure of the
uniform attractor, the continuity maybe   necessary to some
extend.

\subsubsection{Norm-to-weak continuous processes}

In \cite{CV94,CV}, in order to obtain the structure  of the uniform
attractor, the authors assume that the family of processes
$\{U_{\sigma}(t,\tau)\}$, $\sigma\in \Sigma$ is $(X\times \Sigma,
\,X)-$continuous, see Theorem 5.1 in Chapter IV of \cite{CV}.

Now, as noticed in \cite{ZYS}, in order to obtain the invariance of
the global attractor of a semigroup for an autonomous system, we
only need the norm-to-weak continuity. In this part, we will
generalize these results to non-autonomous systems.

\begin{definition}\label{d3.5}
A family of processes $\{U_\sigma (t,\tau)\}$, $\sigma \in \Sigma $
is said to be norm-to-weak continuous, if for any fixed $t$ and
$\tau \in \mathbb{R}$ with $t\geqslant \tau$, for any
$\{x_n\}\subset X$ and $\{\sigma_n\}\subset \Sigma$, we have
\[
\left.
\begin{array}{r}
x_n\overset{\|\cdot\|_X }{\to} x\\
\sigma_n \overset{d}{\to} \sigma
\end{array}
\right\} \Rightarrow U_{\sigma_n}(t,\,\tau)x_n \rightharpoonup
U_{\sigma}(t,\,\tau)x \quad {\rm weakly ~in}~X.
\]
\end{definition}

For convenience, we also use the following notations:
\begin{definition}\label{d3.6}
A semigroup $\{S(t)\}_{t\geqslant 0}$: $X\times \Sigma \to X\times
\Sigma$ is to be called skew productively norm-to-weak continuous,
if for any fixed $t\geqslant 0$, for any $\{x_n\}\subset X$ and
$\{\sigma_n\}\subset \Sigma$, we have
\begin{align*}
&\Pi_1S(t)(x_n,\, \sigma_n)\overset{weak}{\rightharpoonup}
\Pi_1S(t)(x,\,\sigma), \nonumber \\
& \Pi_2S(t)(x_n,\, \sigma_n)\overset{d}{\to} \Pi_2S(t)(x,\,\sigma)
\end{align*}
provided that $x_n\overset{\|\cdot\|_X }{\to} x$ and $\sigma_n
\overset{d}{\to} \sigma$, where $\Pi_1$ and $\Pi_2$ are the
canonical projector from $X\times \Sigma$ to $X$ and $\Sigma$,
respectively. Denoted such continuity by
\[
S(t)(x_n,\, \sigma_n)\overset{s-w}{\rightharpoonup}
S(t)(x,\,\sigma).
\]
\end{definition}

Following from the definition of skew productively norm-to-weak
continuous semigroup, we have the following result.
\begin{proposition}\label{p3.7}
Let $\{U_\sigma (t,\tau)\}$, $\sigma \in \Sigma $ be a family
norm-to-weak continuous processes in $X$, and the translation
semigroup $\{T(t)\}_{t\geqslant 0}$ is continuous (w.r.t. the metric
$d(\cdot,\,\cdot)$) in $\Sigma$. Then, the semigroup
$\{S(t)\}_{t\geqslant 0}$ corresponding to $\{U_\sigma (t,\tau)\}$,
$\sigma \in \Sigma $, defined by \eqref{2.5} and acting on $X\times
\Sigma$, is skew productively norm-to-weak continuous.
\end{proposition}

\subsubsection{Kernel sections  of the uniform attractor}

\begin{theorem}\label{t3.8}
Let $X$ be a Banach space and $\Sigma$ be a compact metric space.
Assume that a family of processes $\{U_\sigma(t,\tau)\}$,
$\sigma\in\Sigma$ satisfies the translation identity
\eqref{2.3}-\eqref{2.4}, as well as the following conditions:
\begin{enumerate}
\item[(i)] The translation semigroup $\{T(t)\}_{t\geqslant 0}$ is
continuous on $\Sigma$;

\item[(ii)] $\{U_\sigma(t,\tau)\},
\sigma\in\Sigma$ is norm-to-weak continuous on $X$;

\item[(iii)] $\{U_\sigma(t,\tau)\},
\sigma\in\Sigma$ has a bounded uniformly (w.r.t. $\sigma \in
\Sigma$) absorbing set $B_0$ in $X$;

\item[(iv)] $\{U_\sigma(t,\tau)\},
\sigma\in\Sigma$ is uniform (w.r.t. $\sigma \in \Sigma$)
asymptotically compact in $X$.
\end{enumerate}
Then, $\{U_\sigma(t,\tau)\}, \sigma\in\Sigma$ has a uniform (w.r.t.
$\sigma \in \Sigma$) attractor $\mathscr{A}_{\Sigma}$ satisfying
\begin{equation}\label{3.3}
\mathscr{A}_{\Sigma}=\omega_{0,\,\Sigma}(B_0)=\bigcup_{\sigma\in
\Sigma}\mathcal{K}_{\sigma}(s), \quad \forall~s\in \mathbb{R},
\end{equation}
where $\mathcal{K}_{\sigma}(s)$ is the section at $t=s$ of the
kernel $\mathcal{K}_{\sigma}$ of the process
$\{U_{\sigma}(t,\,\tau)\}$ with symbol $\sigma$.
\end{theorem}
{\bf Proof.} From the assumptions (iii), (iv) and $Theorem$
\ref{t3.4}, we know that the family of processes
$\{U_\sigma(t,\tau)\},\, \sigma\in\Sigma$ has a compactly uniform
(w.r.t. $\sigma \in \Sigma$) attractor $\mathscr{A}_{\Sigma}$ which
satisfies
\begin{equation*}
\mathscr{A}_{\Sigma}=\omega_{0,\,\Sigma}(B_0).
\end{equation*}

In order to prove the structure \eqref{3.3}, we will construct skew
product flow on $X\times \Sigma$. We will complete the proof by
three steps.

{\it Step 1.} Constructing skew product flow $\{S(t)\}_{t\geqslant
0}$ on $X\times \Sigma$:
\begin{equation}\label{3.4}
S(t)(u,\sigma)=(U_\sigma(t,0)u,T(t)\sigma),\quad t\geqslant
0,\quad(u,\sigma)\in X\times\Sigma.
\end{equation}

Then follows from \eqref{2.3} and $Proposition$ \ref{p3.7} we have
that $\{S(t)\}_{t\geqslant 0}$ forms a skew productively
norm-to-weak continuous semigroup on $X\times \Sigma$.

Define the $\omega-$limit set of $\{S(t)\}_{t\geqslant 0}$ by
\begin{equation*}
\mathscr{A}=\omega(B_0\times \Sigma)=\bigcap_{t\geqslant
0}\overline{\bigcup_{s\geqslant t}S(s)(B_0\times \Sigma)},
\end{equation*}
where $\overline{A}$ denotes the closure of $A$ in $X\times \Sigma$.
Then we have the following equivalent characterization:
\begin{align}\label{3.5}
(x,\,\sigma)\in \mathscr{A}~\text{if and}~ &\text{only if there
exist}\ x_n\in B_0,
\sigma_n\in \Sigma~\text{and}~t_n\geqslant 0~\text{with}~t_n\to \infty, \nonumber \\
&~\text{such that}\ S(t_n)(x_n,\,\sigma_n)\to
(x,\,\sigma)\,\text{as}~n\to \infty~\text{in}~X\times \Sigma.
\end{align}

{\it Step 2.} From the assumptions that $\Sigma$ is a compact metric
space and conditions (iii) and (iv) we have that $\mathscr{A}$ is
nonempty, compact in $X\times \Sigma$, and attracts every bounded
subset of $X\times \Sigma$ under the topology of $X\times \Sigma$.
In order to show that $\mathscr{A}$ is a global attractor of
$\{S(t)\}_{t\geqslant 0}$ in $X\times \Sigma$, we need to show the
invariance of $\mathscr{A}$, that is, we need to prove
\begin{equation*}
S(t)\mathscr{A}=\mathscr{A}\quad \forall~t\geqslant 0.
\end{equation*}

For any $(x,\,\sigma)\in \mathscr{A}$ and $t\geqslant 0$. From
\eqref{3.5} we know that there exist $x_n\in B_0$, $\sigma_n\in
\Sigma$ and $t_n\geqslant 0$, $t_n\to \infty$ such that
$S(t_n)(x_n,\,\sigma_n)\to (x,\,\sigma)$. Since $t_n\to \infty$,
without loss of generality, we assume that $t_n\geqslant t$ for each
$n$ (at most by passing subsequence). Then from \eqref{3.4} we have
\[
S(t_n-t)(x_n,\,\sigma_n)=(U_{\sigma_n}(t_n-t,\,0)x_n,\,T(t_n-t)\sigma_n).
\]

Noticing that $\{x_n\}\subset B_0$ and $t_n-t\to \infty$, then by
the assumption (iv), we know that there is a convergent subsequence
of $\{U_{\sigma_n}(t_n-t,\,0)x_n\}$, without loss of generality, we
also assume that $U_{\sigma_n}(t_n-t,\,0)x_n \to y$ for some $y\in
X$ as $n\to \infty$. At the same time, by the assumption that
$\Sigma$ is a compact metric space, without loss of generality, we
can also assume that $T(t_n-t)\sigma_n \to \sigma_0$ for some
$\sigma_0\in \Sigma$ as $n\to \infty$. Consequently,
$S(t_n-t)(x_n,\,\sigma_n) \to (y,\,\sigma_0)$ as $n\to \infty$. Then
from the definition we have that $(y,\,\sigma_0) \in \mathscr{A}$.
Moreover, from the skew productively norm-to-weak continuity of
$\{S(t)\}_{t\geqslant 0}$, we have
\begin{align*}
(x,\sigma)\leftarrow
S(t_n)(x_n,\,\sigma_n)=S(t)S(t_n-t)(x_n,\,\sigma_n)\overset{s-w}
{\rightharpoonup}S(t)(y,\,\sigma_0).
\end{align*}

Following the uniqueness of limits, we have
$(x,\sigma)=S(t)(y,\,\sigma_0)$, i.e., $\mathscr{A}\subset
S(t)\mathscr{A}$.

Similarly, we can prove $S(t)\mathscr{A}\subset \mathscr{A}$.
Therefore, $\mathscr{A}$ is the global attractor of
$\{S(t)\}_{t\geqslant 0}$ and we have
\begin{equation*}
\mathscr{A}=\{\gamma(0)~|~\gamma(\cdot)~\text{is a bounded complete
trajectory of}~\{S(t)\}_{t\geqslant 0}\}.
\end{equation*}

{\it Step 3.} Based on Step 1 and Step 2, similar to what done in
Chepyzhov \& Vishik [\cite{CV}, $Theorem$ 5.1, Chapter IV], we can
prove that $\mathscr{A}_{\Sigma}=\Pi_1\mathscr{A}=\bigcup_{\sigma\in
\Sigma}\mathcal{K}_{\sigma}(s)$ and $\Pi_2\mathscr{A}=\Sigma$.
$\hfill \blacksquare$

In practical applications, $\Sigma$ would be the completion of a
dense subset $\Sigma_0(\subset \Sigma)$ with respect to some metric
$d(\cdot,\,\cdot)$, and maybe different for different metrics. For
example, $\Sigma_0=\{g_0(s+h,x)|~h\in \mathbb{R}\}$ for some
function $g_0(s,x)$ belong to a special function space, and $\Sigma$
can be chosen according to our concrete problem.

If the family of processes $\{U_{\sigma}(t,\,\tau)\},\sigma\in
\Sigma_0$ satisfies also the translation identity
\eqref{2.3}-\eqref{2.4}, then, Chepyzhov \& Vishik proved in
\cite{CV94,CV} that the uniform (w.r.t $\tau\in \mathbb{R}$)
attractor of a process $\{U_{\sigma_0}(t,\tau)\}$, $\tau \in
\mathbb{R}$ coincides with the uniform (w.r.t. $\sigma\in \Sigma_0$)
attractor for the family of processes
$\{U_{\sigma}(t,\,\tau)\},\sigma\in \Sigma_0$.

Similar to \cite{CV,LWZ}, the following results give a method to
obtain the structure of the uniform (w.r.t $\tau\in \mathbb{R}$)
attractor of a process $\{U_{\sigma_0}(t,\tau)\}$, $\tau \in
\mathbb{R}$ via the structure of the uniform (w.r.t. $\sigma\in
\Sigma_0$) attractor for the family of processes
$\{U_{\sigma}(t,\,\tau)\},\sigma\in \Sigma_0$.

Since our processes are norm-to-weak continuity, we first give a
simple lemma about metrizable. We recall (e.g., see
Diestel\cite{Diestel}, p. 18) that a set $F\subset X^*$ is called
total if $f(x)=0$ for each $f \in F$ implies $x=0$.
\begin{lemma}\label{a}
If $K$ is a (relatively) weakly compact subset in a Banach space $X$
and $K$ is countable, then $\overline{K}^{weak}$ is metrizable,
where $\overline{K}^{weak}$ means the weak closure of $K$ in $X$.
\end{lemma}
{\bf Proof.} Denote $Y=\overline{span\{K\}}$.

From the convexity of $span\{K\}$ we know $Y$ is weakly closed in
$X$. Therefore, $K=K\cap Y$ is (relatively) weakly compact in the
separable Banach space $Y$. Since the dual of a separable Banach
space contains a countable total set, we know that
$\overline{K}^{weak}$ is metrizable in $Y$, and from that $Y$ is a
closed subspace of $X$ we get $\overline{K}^{weak}$ is metrizable in
$X$. $\hfill \blacksquare$

\begin{theorem}\label{t3.9}
Let $\Sigma_0$ be a parameter set, $\Sigma$ is a completion of
$\Sigma_0$ with respect to some metric $d(\cdot,\,\cdot)$, and the
translation semigroup $\{T(t)\}_{t\geqslant 0}$ satisfies also the
translation identity \eqref{2.3}-\eqref{2.4} on $\Sigma_0$.
Furthermore, assume that the family of processes
$\{U_{\sigma}(t,\,\tau)\},\sigma\in \Sigma$ satisfies all of the
assumptions in Theorem \ref{t3.8}. Then, both families of processes
$\{U_{\sigma}(t,\,\tau)\}$,$\sigma\in \Sigma$ and $\sigma\in
\Sigma_0$ have compactly uniform (w.r.t. $\sigma\in \Sigma$ and
$\sigma\in \Sigma_0$ respectively) attractors $\mathscr{A}_{\Sigma}$
and $\mathscr{A}_{\Sigma_0}$ respectively, moreover,
\begin{equation*}
\mathscr{A}_{\Sigma_0}=\mathscr{A}_{\Sigma}=\omega_{0,\,\Sigma}(B_0)=\bigcup_{\sigma\in
\Sigma}\mathcal{K}_{\sigma}(s), \quad \forall~s\in \mathbb{R}.
\end{equation*}
\end{theorem}
{\bf Proof.} The existence is a immediate  consequence of
$Theorem$ \ref{t3.8}, and obviously, we have
\begin{equation*}
\mathscr{A}_{\Sigma_0}\subset\mathscr{A}_{\Sigma}=\omega_{0,\,\Sigma}(B_0)=\bigcup_{\sigma\in
\Sigma}\mathcal{K}_{\sigma}(s), \quad \forall~s\in \mathbb{R}.
\end{equation*}

Now we prove $\omega_{0,\Sigma_0}(B_0)=\omega_{0,\Sigma}(B_0)$. For
any $y\in \omega_{0,\,\Sigma}(B_0)$, from $Proposition$ \ref{p3.2},
we know that there exist $x_n\in B_0$, $t_n\to\infty$ and
$\sigma_n\in \Sigma$ such that
\begin{equation}\label{3.6}
U_{\sigma_n}(t_n,0)x_n\to y\quad \text{as}~n\to \infty.
\end{equation}
On the other hand, from the assumption that $\Sigma$ is the
completion of $\Sigma_0$ we know that there exists
$\{\sigma^{(n)}_m\}\subset \Sigma_0$ satisfies
$\sigma^{(n)}_m\overset{d}{\to}\sigma_n$ as $m\to \infty$ for each
$n\in\mathbb{N}$. Therefore, due to the norm-to-weak continuity of
the family of processes $\{U_{\sigma}(t,\,\tau)\}$, $\sigma\in
\Sigma$, we have
\begin{equation}\label{3.7}
U_{\sigma^{(n)}_m}(t_n,0)x_n \rightharpoonup
U_{\sigma_n}(t_n,0)x_n\quad \text{as} ~m\to \infty
\end{equation}
for each $n\in \mathbb{N}$. Denote
$K=\{U_{\sigma^{(n)}_m}(t_n,0)x_n~|n,m\in \mathbb{N}\}$, then $K$ is
countable and thanks to the condition (iv) of $Theorem$ \ref{t3.8}
we know $K$ is also relatively weakly compact in $X$. Consequently,
from $Lemma$ \ref{a} we have that $\overline{K}^{weak}$ is
metrizable.

Hence, combining \eqref{3.6} and \eqref{3.7}, we can obtain that
there exist $\sigma'_n\in \Sigma_0$ for each $n\in \mathbb{N}$ such
that
\begin{equation*}
U_{\sigma'_n}(t_n,0)x_n \rightharpoonup y\quad
\text{in}~X~\text{as}~n\to \infty.
\end{equation*}
Then noticing the uniform asymptotic compactness again and the
uniqueness of limits, we have $y\in \omega_{0,\Sigma_0}(B_0)$.
$\hfill \blacksquare$

\section{A criterion for verifying the uniform asymptotic compactness}

In this section, we present a technical method to verify the uniform
asymptotic compactness (given in $Definition$ 3.1) for the family of
processes generated by non-autonomous hyperbolic type of
evolutionary equations. This criterion is partially motivated by the
methods in \cite{CL1,CL2,CL3,Kh,SYZ2} for autonomous systems.
Here, the following results
and proof are similar to that in \cite{SYZ2} for autonomous cases.

\begin{definition}
Let $X$ be a Banach space and $B$ be a bounded subset of $X$,
$\Sigma$ be a symbol (or parameter) space. We call a function
$\phi(\cdot,\cdot;\,\cdot,\cdot)$, defined on $(X\times
X)\times(\Sigma \times\Sigma)$, to be a contractive function on
$B\times B$ if for any sequence $\{x_n\}_{n=1}^{\infty}\subset B$
and any $\{\sigma_n\}\subset \Sigma$, there is a subsequence
$\{x_{n_k}\}_{k=1}^{\infty} \subset \{x_n\}_{n=1}^{\infty}$ and
$\{\sigma_{n_k}\}_{k=1}^{\infty} \subset
\{\sigma_n\}_{n=1}^{\infty}$ such that
\begin{equation*}
\lim_{k\to \infty}\lim_{l\to
\infty}\phi(x_{n_k},\,x_{n_l};\,\sigma_{n_k},\,\sigma_{n_l})=0.
\end{equation*}
We denote the set of all contractive functions on $B\times B$ by
$\mathfrak{C}(B,\,\Sigma)$.
\end{definition}

\begin{theorem}\label{t4.2}
Let $\{U_{\sigma}(t,\tau)\}$, $\sigma\in\Sigma$ be a family of
processes satisfies the translation identity \eqref{2.3}-\eqref{2.4}
on Banach space $X$ and has a bounded uniformly (w.r.t. $\sigma \in
\Sigma$) absorbing set $B_0\subset X$. Moreover, assume that for any
$\varepsilon >0$ there exist $T=T(B_{0}, \varepsilon)$ and
$\phi_T\in\mathfrak{C}(B_0,\,\Sigma)$ such that
\begin{equation*}
\|U_{\sigma_1}(T,0)x-U_{\sigma_2}(T,0)y\|\leqslant\varepsilon+\phi_T(x,y;\sigma_1,\sigma_2),\quad
\forall\ \,x,\,y\in B_0,~\forall\ \sigma_1,\,\sigma_2\in \Sigma.
\end{equation*}
Then $\{U_\sigma(t,\tau)\}$, $\sigma\in\Sigma$ is uniformly (w.r.t.
$\sigma\in \Sigma$) asymptotically compact in $X$.
\end{theorem}
{\bf Proof.} For any fixed $\tau \in\mathbb{R}$, let
$\{x_n\}_{n=1}^{\infty}$ be a bounded sequence of $X$, $\sigma_n\in
\Sigma$ and $t_n \geqslant \tau$ satisfy $t_n\to \infty$ as $n\to
\infty$. We need to show that
\begin{align*}
\{U_{\sigma_n}(t_n,\tau )x_n\}_{n=1}^\infty~\text{is precompact in}~
X.
\end{align*}

Thanks to the translation identity \eqref{2.3}-\eqref{2.4}, we know
that for any fixed $\tau\in\mathbb{R}$ and $\sigma\in \Sigma$ we can
find $\sigma'\in\Sigma$ such that
\begin{equation}\label{4.1}
U_{\sigma'}(t+\tau,\tau)x=U_\sigma(t,0)x,\qquad \text{for
all}~t\geqslant 0 ~\text{and}~x\in X.
\end{equation}
Therefore, we only need to show that
$\{U_{\sigma_n}(t_n,0)x_n\}_{n=1}^{\infty}$ is precompact in $X$.

In the following, we will prove that
$\{U_{\sigma_n}(t_n,0)x_n\}_{n=1}^{\infty}$ has a Cauchy subsequence
via a diagonal method.

Taking $\varepsilon_m>0$ with $\varepsilon_m \to 0$ as $m\to
\infty$.

At first, for $\varepsilon_1$, by the assumptions, there exist
$T_1=T_1(\varepsilon_1)$ and $\phi_1 \in \mathfrak{C}(B_0,\Sigma)$
such that
\begin{equation}\label{4.2}
\|U_{\sigma_1}(T_1,0)x-U_{\sigma_2}(T_1,0)y\|\leqslant\varepsilon_1+\phi_1(x,y;\sigma_1,
\sigma_2)~ \text{for all}~x,y\in B_0~ \text{and}~
\sigma_1,\sigma_2\in \Sigma.
\end{equation}

Since $t_n\to \infty$, for such fixed $T_1$, without loss of
generality, we assume that $t_n\gg T_1$ is so large that
$U_{\sigma_n}(t_n-T_1,0)x_n \in B_0$ for each $n\in \mathbb{N}$.

Similar to \eqref{4.1}, for each $n\in \mathbb{N}$, there is a
$\sigma'_n\in \Sigma$ such that
\begin{equation}\label{4.3}
U_{\sigma'_n}(T_1,0)=U_{\sigma_n}(t_n,t_n-T_1).
\end{equation}

Let $y_n=U_{\sigma_n}(t_n-T_1,0)x_n$, then from \eqref{4.2} and
\eqref{4.3} we have
\begin{align}\label{4.4}
\|U_{\sigma_n}&(t_n,0)x_n-U_{\sigma_m}(t_m,0)x_m\| \nonumber \\ & =
\|U_{\sigma_n}(t_n,t_n-T_1)U_{\sigma_n}(t_n-T_1,0)x_n-U_{\sigma_m}(t_m,t_m-T_1)
U_{\sigma_n}(t_m-T_1,0)x_m\|
\nonumber \\
&= \|U_{\sigma_n}(t_n,t_n-T_1)y_n-U_{\sigma_m}(t_m,t_m-T_1)y_m\|
\nonumber \\
&=\|U_{\sigma'_n}(T_1,0)y_n-U_{\sigma'_m}(T_1,0)y_m\| \nonumber \\
&\leqslant\varepsilon_1+\phi_1(y_n,y_m;\sigma'_n,\sigma'_m).
\end{align}

Due to the definition of $\mathfrak{C}(B_0,\Sigma)$ and $\phi_1 \in
\mathfrak{C}(B_0,\Sigma)$, we know that $\{y_n\}_{n=1}^{\infty}$ has
a subsequence $\{y^{(1)}_{n_k}\}_{k=1}^{\infty}$ and
$\{\sigma'_n\}_{n=1}^{\infty}$ has a subsequence
$\{\sigma'^{(1)}_{n_k}\}_{k=1}^{\infty}$ such that
\begin{equation}\label{4.5}
\lim_{k\to \infty}\lim_{l\to
\infty}\phi_1(y^{(1)}_{n_k},\,y^{(1)}_{n_l};\sigma'^{(1)}_{n_k},\sigma'^{(1)}_{n_l})
\leqslant \frac{\varepsilon_1}{2}.
\end{equation}
And similar to the autonomous cases, e.g., see \cite{Kh,SYZ2}, we
have
\begin{align*}
\lim_{k\to \infty}\sup_{p\in \mathbb{N}}&
\|U_{\sigma^{(1)}_{n_{k+p}}}(t^{(1)}_{n_{k+p}},\,0)x^{(1)}_{n_{k+p}}
-U_{\sigma^{(1)}_{n_{k}}}(t^{(1)}_{n_{k}},\,0)x^{(1)}_{n_{k}}\|
\nonumber \\
&\leqslant \lim_{k\to \infty}\sup_{p\in \mathbb{N}}\limsup_{l\to
\infty}
\|U_{\sigma^{(1)}_{n_{k+p}}}(t^{(1)}_{n_{k+p}},\,0)x^{(1)}_{n_{k+p}}
-U_{\sigma^{(1)}_{n_{l}}}(t^{(1)}_{n_{l}},\,0)x^{(1)}_{n_{l}}\|
\nonumber \\
& \qquad + \limsup_{k\to \infty}\limsup_{l\to \infty}
\|U_{\sigma^{(1)}_{n_{k}}}(t^{(1)}_{n_{k}},\,0)x^{(1)}_{n_{k}}
-U_{\sigma^{(1)}_{n_{l}}}(t^{(1)}_{n_{l}},\,0)x^{(1)}_{n_{l}}\|
\nonumber \\
&\leqslant \varepsilon_1+\lim_{k\to \infty}\sup_{p\in
\mathbb{N}}\lim_{l\to
\infty}\phi_1(y^{(1)}_{n_{k+p}},\,y^{(1)}_{n_l};\sigma'^{(1)}_{n_{k+p}},\sigma'^{(1)}_{n_l})
\nonumber \\
& \qquad +\varepsilon_1+ \lim_{k\to \infty}\lim_{l\to
\infty}\phi_1(y^{(1)}_{n_k},\,y^{(1)}_{n_l};\sigma'^{(1)}_{n_k},\sigma'^{(1)}_{n_l}),
\end{align*}
which, combining with \eqref{4.4} and \eqref{4.5}, implies that
\begin{align*}
\lim_{k\to \infty}\sup_{p\in \mathbb{N}}&
\|U_{\sigma^{(1)}_{n_{k+p}}}(t^{(1)}_{n_{k+p}},\,0)x^{(1)}_{n_{k+p}}
-U_{\sigma^{(1)}_{n_{k}}}(t^{(1)}_{n_{k}},\,0)x^{(1)}_{n_{k}}\|
\leqslant 4\varepsilon_1.
\end{align*}

Therefore, there is a $K_1$ such that
\begin{align*}
\|U_{\sigma^{(1)}_{n_{k}}}(t^{(1)}_{n_{k}},\,0)x^{(1)}_{n_{k}}
-U_{\sigma^{(1)}_{n_{l}}}(t^{(1)}_{n_{l}},\,0)x^{(1)}_{n_{l}}\|
\leqslant 5\varepsilon_1 \quad \forall~k,l\geqslant K_1.
\end{align*}

By induction, we obtain that, for each $m\geqslant 1$, there is a
subsequence $\{U_{\sigma^{(m+1)}_{n_{k}}}
(t^{(m+1)}_{n_{k}},\,0)x^{(m+1)}_{n_{k}}\}_{k=1}^{\infty}$ of
$\{U_{\sigma^{(m)}_{n_{k}}}(t^{(m)}_{n_{k}},\,0)x^{(m)}_{n_{k}}\}_{k=1}^{\infty}$
and certain $K_{m+1}$ such that
\begin{equation*}
\|U_{\sigma^{(m+1)}_{n_{k}}}(t^{(m+1)}_{n_{k}},\,0)x^{(m+1)}_{n_{k}}
-U_{\sigma^{(m+1)}_{n_{l}}}(t^{(m+1)}_{n_{l}},\,0)x^{(m+1)}_{n_{l}}\|
\leqslant 5\varepsilon_{m+1} \quad \forall~k,l\geqslant K_{m+1}.
\end{equation*}

Now, we consider the diagonal subsequence
$\{U_{\sigma^{(k)}_{n_{k}}}(t^{(k)}_{n_{k}},\,0)x^{(k)}_{n_{k}}\}_{k=1}^{\infty}$.
Since for each $m\in \mathbb{N}$,
$\{U_{\sigma^{(k)}_{n_{k}}}(t^{(k)}_{n_{k}},\,0)x^{(k)}_{n_{k}}\}_{k=m}^{\infty}$
is a subsequence of
$\{U_{\sigma^{(m)}_{n_{k}}}(t^{(m)}_{n_{k}},\,0)x^{(m)}_{n_{k}}\}_{k=1}^{\infty}$,
then,
\begin{equation*}
\|U_{\sigma^{(k)}_{n_{k}}}(t^{(k)}_{n_{k}},\,0)x^{(k)}_{n_{k}}
-U_{\sigma^{(l)}_{n_{l}}}(t^{(l)}_{n_{l}},\,0)x^{(l)}_{n_{l}}\|
\leqslant 5\varepsilon_{m} \quad \forall~k,l\geqslant
\max\{m,K_{m}\},
\end{equation*}
which, combining with $\varepsilon_m \to 0$ as $m\to \infty$,
implies that
$\{U_{\sigma^{(k)}_{n_{k}}}(t^{(k)}_{n_{k}},\,0)x^{(k)}_{n_{k}}\}_{k=m}^{\infty}$
is a Cauchy sequence in $X$. This shows that
$\{U_{\sigma_n}(t_n,0)x_n\}_{n=1}^{\infty}$ is precompact in $X$.
$\hfill \blacksquare$

\section{Application to wave equation}

\subsection{Mathematical setting}

Similar to the autonomous cases (e.g., see \cite{CL2}), applying the
Galerkin approximation method, we have the following existence and
uniqueness results (e.g., see \cite{Ha2,Li}), and the time-dependent
terms make no essential complications.
\begin{theorem}\label{t5.1}
Let $\Omega$ be a bounded domain of $\mathbb{R}^3$ with smooth
boundary, $h$ and $f$ satisfy \eqref{1.4}-\eqref{1.8}, and $g\in
L^{\infty}(\mathbb{R};\,L^2(\Omega))$. Then the non-autonomous
initial-boundary value problem \eqref{1.1}-\eqref{1.3} has an unique
solution $u(t)$ satisfying $(u(t),\,u_t(t))\in
\mathcal{C}(\mathbb{R}^{\tau};\,H_0^1(\Omega)\times L^2(\Omega))$
and $\partial_{tt}u(t)\in
L^2_{loc}(\mathbb{R}^{\tau};\,H^{-1}(\Omega))$ for any initial data
$(u^{0\tau},\,u^{1\tau})\in H_0^1(\Omega)\times L^2(\Omega)$.
\end{theorem}

We use the notations as in Chepyzhov \& Vishik\cite{CV}: Let
$y(t)=(u(t),\,u_t(t))$, $y_{\tau}=(u^{0\tau},\,u^{1\tau})$ and
$X=H_0^1(\Omega)\times L^2(\Omega)$ with finite energy norm
\begin{equation*}
\|y\|_{X}=\{\|\nabla u\|^2+|u_t|^2\}^{\frac{1}{2}}.
\end{equation*}

Let $A_{\sigma(t)}(u,v)=(v,\Delta u-f(u)-h(v)+\sigma(t))$. Then
the non-autonomous system \eqref{1.1}-\eqref{1.3} can be rewritten
in the operator form
\begin{equation}\label{5.1}
\partial_t y=A_{\sigma(t)}(y),\qquad y|_{t=\tau}=y_\tau,
\end{equation}
where $\sigma(s)=g(x,s)$ is symbol of equation \eqref{5.1}.

We now define the symbol space for \eqref{5.1}. Taking a fixed
symbol $\sigma_0(s)=g_0(x,s)$, $g_0 \in
L^{\infty}(\mathbb{R};\,L^2(\Omega))\cap
W_{b}^{1,r}(\mathbb{R};\,L^r(\Omega))$. Set
\begin{equation}\label{5.2}
\Sigma_0=\{g_0(x,\,t+h)~|~h\in \mathbb{R}\}
\end{equation}
and
\begin{equation}\label{5.3}
\Sigma~\text{be the $*-$weakly closure of}~\Sigma_0 ~\text{in}~
L^{\infty}(\mathbb{R};\,L^2(\Omega)) \cap
W_{b}^{1,\,r}(\mathbb{R};\,L^r(\Omega)).
\end{equation}

Then we have the following simple properties.
\begin{proposition}\label{p5.2}
\end{proposition}\vspace{-0.5 cm}
\begin{enumerate}
\item[(i)]\ $\Sigma$ is bounded in
$L^{\infty}(\mathbb{R};\,L^2(\Omega)) \cap
W_{b}^{1,\,r}(\mathbb{R};\,L^r(\Omega))$, and for any $\sigma\in
\Sigma$, the following estimate holds
\begin{equation*}
\|\sigma\|_{L^{\infty}(\mathbb{R};\,L^2(\Omega)) \cap
W_{b}^{1,\,r}(\mathbb{R};\,L^r(\Omega))} \leqslant
\|g_0\|_{L^{\infty}(\mathbb{R};\,L^2(\Omega)) \cap
W_{b}^{1,\,r}(\mathbb{R};\,L^r(\Omega))};
\end{equation*}
\item[(ii)] The translation semigroup $\{T(h)| h\geqslant 0\}$
acting on $\Sigma$ is invariant in $\Sigma$, that is
\begin{equation*}
T(h)\Sigma=\Sigma\quad \text{for all} ~h\in\mathbb{R}^+.
\end{equation*}
\end{enumerate}

Thus, from $Theorem$ \ref{t5.1}, we know that
\eqref{1.1}-\eqref{1.3} is well posed for all $\sigma(s)\in
\Sigma$ and generates a family of processes
$\{U_{\sigma}(t,\tau)\},\sigma\in \Sigma$ given by the formula
$U_{\sigma}(t,\tau)y^{\tau}=y(t)$, where $y(t)$ is the solution of
\eqref{1.1}-\eqref{1.8}, and $\{U_{\sigma}(t,\tau)\},\sigma\in
\Sigma$ satisfies \eqref{2.1}-\eqref{2.2}. At the same time, by
the unique solvability, we know $\{U_{\sigma}(t,\tau)\},\sigma\in
\Sigma$ satisfies the translation identity \eqref{2.3}.

In what follows, we denote by $\{U_{\sigma}(t,\tau)\},\sigma\in
\Sigma$ the family of processes generated by
\eqref{5.1}-\eqref{5.3}.

\subsection{Bounded uniformly (w.r.t. $\sigma\in \Sigma$) absorbing set}
We begin with the following result on the existence of bounded
uniformly (w.r.t. $\sigma\in \Sigma$) absorbing set. Its   proof
is essentially established in Haraux\cite{Ha2}, and for reader's
convenience, we replicate it here and only make a few minor
changes for our problem.

\begin{theorem}\label{t5.5}
Under the assumptions of Theorem \ref{t5.1}, the family of
processes $\{U_{\sigma}(t,\tau)\}$, $\sigma\in \Sigma$
corresponding to \eqref{5.1} has a bounded (in $X$) uniformly
(w.r.t. $\sigma\in \Sigma$) absorbing set $B_0$, i.e., there
exists a positive constant $\rho$, which depends on
$\|g_0\|_{L^{\infty}(\mathbb{R};\,L^2(\Omega))}$ and the
coefficients in \eqref{1.6}-\eqref{1.8}, such that for any bounded
subset $B\subset X$ and any $\tau \in \mathbb{R}$, there is a
$T=T(B)$ such that for any $t-\tau\geqslant T$, $\sigma\in \Sigma$
and $(u^{0\tau},\,u^{1\tau})\in B$,
\begin{equation*}
\|U_{\sigma}(t,\tau)(u^{0\tau},\,u^{1\tau})\|_X\leqslant \rho.
\end{equation*}
\end{theorem}
{\bf Proof.} Since $\{U_{\sigma}(t,\tau)\}$, $\sigma\in \Sigma$
satisfies the translation identity, we only need to prove $Theorem$
\ref{t5.5} for the cases $\tau\equiv 0$. Moreover, from the
definition of $\Sigma$ we know that for all $\sigma\in \Sigma$,
\[
\|\sigma\|_{L^{\infty}(\mathbb{R};L^2(\Omega))}\leqslant
\|g_0\|_{L^{\infty}(\mathbb{R};L^2(\Omega))}.
\]
Hence, without loss of generality, in the remainder of proof, we
will not point out the difference of symbols and will denote
different $\sigma$ by $g$.

For any $\varepsilon\geqslant0$, we set
\begin{equation}\label{5.4}
E_\varepsilon(t)=\frac12\|u(t)\|^2+\frac12|u_t(t)|^2+\int_\Omega
F(u(x))dx+\varepsilon\langle u_t(t),u(t)\rangle.
\end{equation}
Then we have $E_\varepsilon(t)\rightarrow E_0(t)$ as $\varepsilon\to
0$. Moreover there exist $C_0,\, C_1\geqslant0$ such that
\begin{equation}\label{5.5}
\frac{C_0}{2}(\|u(t)\|^2+|u_t(t)|^2)-C_1\leqslant E_0(t).
\end{equation}
By differentiating \eqref{5.4} with time $t$, we obtain that
\begin{align*}
\dfrac{d}{dt}(E_\varepsilon(t)) =&\langle u_{tt}-\bigtriangleup
u,u_t\rangle+\langle
f(u),u_t(t)\rangle+\varepsilon|u_t(t)|^2+\varepsilon\langle
u_{tt}(t),u_t(t)\rangle\\
=&\langle g,u_t\rangle-\langle
h(u_t),u_t\rangle+\varepsilon|u_t(t)|^2-\varepsilon\|u\|^2-\varepsilon\langle
f(u),u\rangle+\varepsilon\langle g,u\rangle-\varepsilon\langle
h(u_t),u\rangle.
\end{align*}

It is obviously that \eqref{1.5} implies that
\begin{equation*}
\langle h(u_t),u_t\rangle\geqslant\alpha|u_t|^2-C_2|\Omega|,
\end{equation*}
and from \eqref{1.7} and \eqref{1.8} we know that there are
$\lambda_1>\delta>0$ and $C>0$ such that
\begin{equation*}
\langle f(u),u\rangle\geqslant\delta\int_\Omega F(u)dx-C.
\end{equation*}
Hence we get the following inequality
\begin{align*}
E'_\varepsilon(t)\leqslant&(2\varepsilon-\frac
{\alpha}{2})|u_t(t)|^2-\frac{\varepsilon}{2}\|u(t)\|^2-\varepsilon\delta\int_\Omega
F(u)dx\\
&-\frac12\langle
h(u_t),u_t\rangle+\varepsilon\|u(t)\|\|h(u_t(t))\|_{H^{-1}}+C,
\end{align*}
where $C$ depends on $\|g\|_{L^\infty(\mathbb{R},L^2(\Omega))}$.

On the other hand, from \eqref{1.4}-\eqref{1.6} we have (e.g., see
the $Lemma$ in \cite{Ha2}) that there is a constant $K$ such that
\begin{equation*}
\|h(v)\|_{H^{-1}}\leqslant K(1+\langle h(v),v\rangle) \quad
\text{for all}~v\in H_0^1(\Omega).
\end{equation*}

Denote
\begin{equation*}
w(t)\triangleq 1+\langle h(u_t(t)),u_t(t)\rangle(>0).
\end{equation*}
Then, by taking $\varepsilon$ small enough, we obtain that, for all
$t\geqslant0$,
\begin{align}\label{5.6}
E'_\varepsilon(t)&\leqslant-\gamma\varepsilon
E_\varepsilon(t)+(K\varepsilon\|u(t)\|-1/2)w(t)+C \nonumber \\
&\leqslant-\gamma\varepsilon
E_\varepsilon(t)+(N\varepsilon\sqrt{E_\varepsilon(t)}-1/2)w(t)+C,
\end{align}
where $N,C>0$ depending only on $f,g,h$ and $\Omega$ (not on the
initial data) and $\gamma>0$.

Now choose $\varepsilon>0$ so small that
\begin{equation*}
E_\varepsilon(0)<(\frac{1}{2N\varepsilon})^2-\frac{C}{\gamma\varepsilon}.
\end{equation*}
Then,
\begin{equation}\label{5.7}
E_\varepsilon(t)<(\frac{1}{2N\varepsilon})^2,\qquad\forall\,
t\geqslant 0.
\end{equation}

If \eqref{5.7} is not true, let $t_0=\inf\{t\geqslant0,
E_\varepsilon(t)\geqslant(\frac{1}{2N\varepsilon})^2\}$, then
$E_\varepsilon(t_0)=(\frac{1}{2N\varepsilon})^2$ and for all
$t\in[0,t_0]$, we have
\begin{equation}\label{5.8}
E_\varepsilon(t)\leqslant(\frac{1}{2N\varepsilon})^2.
\end{equation}
Therefore, from \eqref{5.6} and \eqref{5.8} we can obtain that
\begin{equation*}
E_\varepsilon(t_0)\leqslant e^{-\gamma\varepsilon t_0}
E_\varepsilon(0)+\frac{C}{\gamma\varepsilon}<
(\frac{1}{2N\varepsilon})^2.
\end{equation*}
This is a contradiction and means that \eqref{5.7} is indeed
satisfies.

Combining \eqref{5.6} and \eqref{5.7}, by use of the uniform
Gronwall lemma, we obtain that
\begin{equation*}
E_\varepsilon(t)\leqslant e^{-\gamma\varepsilon t}
E_\varepsilon(0)+\frac{C}{\gamma\varepsilon}.
\end{equation*}

Finally, we notice that for every bounded set $B\subset X$, assume
the bounds of $B$ (in $X$) is $E(>0)$, then by taking
$1/\varepsilon=4N\sqrt E$ we can obtain
$\frac{1}{(2N\varepsilon)^2}-\frac{C}{\gamma\varepsilon}\geqslant E$
for any $E$ large enough. It follows that there exist $M>0$
(independent of the initial data) and $T=T(B)$ such that
\begin{equation}\label{5.9}
E_0(t)\leqslant M(1+\sqrt{E_0(0)}) \quad \text{for all}~t\geqslant
T~\text{and}~(u^{0\tau},u^{1\tau})\in B.
\end{equation}

Without loss of generality, assume $M>1$, then from \eqref{5.9} we
have that for any bounded set $B\subset X$, there is a $T=T(B)$ such
that
\begin{equation*}
E_0(t)\leqslant 4M^2+1 \quad \text{for all}~t\geqslant
T~\text{and}~(u^{0\tau},u^{1\tau})\in B.
\end{equation*}
Combining with \eqref{5.5} we know that $Theorem$ \ref{t5.5} is
true. \hfill$\blacksquare$

\subsection{Uniform (w.r.t. $\sigma\in \Sigma$) asymptotic
compactness}

The main result in this subsection is summarized in the following
theorem.
\begin{theorem}\label{t5.6}
Let $\Omega$ be a bounded domain in $\mathbb{R}^3$ with smooth
boundary, and $h$ and $f$ satisfy \eqref{1.4}-\eqref{1.8}. If
$g_0\in L^{\infty}(\mathbb{R};\,L^2(\Omega)) \cap
W_{b}^{1,\,r}(\mathbb{R};\,L^r(\Omega))$ and $\Sigma$ is defined
by \eqref{5.3}, then the family of processes
$\{U_{\sigma}(t,\tau)\}$, $\sigma\in \Sigma$ corresponding to
\eqref{5.1} or \eqref{1.1}, is uniformly (w.r.t. $\sigma\in
\Sigma$) asymptotically compact in $H_0^1(\Omega)\times
L^2(\Omega)$.
\end{theorem}

Hereafter, we always assume that the hypotheses of $Theorem$
\ref{t5.1} hold and denote by $B_0$ the bounded uniformly absorbing
set obtained in $Theorem$ \ref{t5.5}.

\subsubsection{Preliminaries}

Note that condition \eqref{1.6} implies that
\[
|h(s)|^{\frac{1}{p}}\leqslant C(1+|s|),
\]
therefore, we have
\[
|h(s)|^{\frac{p+1}{p}}=|h(s)|^{\frac{1}{p}}\cdot |h(s)|\leqslant
C(1+|s|)|h(s)|\leqslant C|h(s)|+Ch(s)\cdot s.
\]
Combining Young inequality and \eqref{1.4} we obtain that
\begin{equation}\label{5.10}
|h(s)|^{\frac{p+1}{p}}\leqslant C(1+h(s)\cdot s)\quad \text{for
all}~s\in \mathbb{R},
\end{equation}
where the constant $C$ is independent of $s$.  And from \eqref{1.4}
and \eqref{1.5}, we have also that
\begin{equation}\label{5.11}
\frac {C 1}{2}s^2\leqslant h(s)s+C \quad \text{for all}~s\in
\mathbb{R}.
\end{equation}
Moreover, we recall the following result.
\begin{lemma}\cite{Fe,Kh}\label{1}
Let $h$ satisfy \eqref{1.4} and \eqref{1.5}. Then for any $\delta
>0$, there exists a constant $C_{\delta}$,   depending on
$\delta$,
such that
\[
|u-v|^2\leqslant \delta + C_{\delta}(h(u)-h(v))(u-v)\quad \text{for
any}~u,v\in \mathbb{R}.
\]
\end{lemma}

\begin{proposition}\label{p5.8}
Let $s_i\in \mathbb{R}$ $(i=1,2,\cdots)$ and $g \in
L^{\infty}(\mathbb{R};\,L^2(\Omega))$. Then there exists $M>0$
such that
\begin{equation*}
\|g(x,\,s_i+t)\|_{L^2(\Omega)} \leqslant M \quad \text{for all}~t\in
\mathbb{R}\setminus \Lambda~\text{and all}~i=1,2,\cdots,
\end{equation*}
where $\Lambda \subset \mathbb{R}$ with $mes(\Lambda)=0$ in
$\mathbb{R}$.
\end{proposition}
{\bf Proof.} Since $g\in L^{\infty}(\mathbb{R};\,L^2(\Omega))$, we
know that there is an $M>0$ such that for each $s_i$,
\begin{equation*}
\|g(x,\,t+s_i)\|_{L^2(\Omega)} \leqslant M\quad \text{for all}~t\in
\mathbb{R}\setminus \Lambda_i,
\end{equation*}
where $mes(\Lambda_i)=0$ in $\mathbb{R}$.

Then $Proposition$ \ref{p5.8} follows immediately by taking
$\Lambda=\bigcup_{i=1}^{\infty}\Lambda_i$. $\hfill \blacksquare$\\

Applying $Proposition$ 7.1 of Robinson\cite{Ro} and $Proposition$
\ref{p5.8} above we can deduce the following results immediately:
\begin{proposition}\label{p5.9}
Let $g \in L^{\infty}(\mathbb{R};\,L^2(\Omega))\cap
W_{b}^{1,r}(\mathbb{R};\,L^r(\Omega))(r>\frac 65)$. Then there is an
$M>0$ such that
\begin{align*}
\sup_{t\in \mathbb{R}}\|g(x,\,t+s)\|_{L^2(\Omega)} \leqslant M\quad
\text{for all}~s\in\mathbb{R}.
\end{align*}
\end{proposition}

\begin{proposition}\label{p5.10}
Let $s_i\in \mathbb{R}$ $(i=1,2,\cdots)$ and $g \in
L^{\infty}(\mathbb{R};\,L^2(\Omega))\cap
W_{b}^{1,r}(\mathbb{R};\,L^r(\Omega))(r>\frac 65)$. Then there
exists $M>0$ such that for any $w\in
W_{loc}^{1,\,2}(\mathbb{R};\,L^2(\Omega))$ and any $T>0$,
\begin{align*}
|\int_0^T\int_s^T&\int_{\Omega}(g(x,\tau+s_i)-g(x,\tau+s_j))w_t(\tau)dxd\tau
ds| \nonumber \\
&\leqslant 2TM\|w(T)\|_{L^2(\Omega)}+2MT^{\frac
12}\left(\int_0^T\int_{\Omega}|w(s)|^2dxds\right)^{\frac 12}
\nonumber
\\
&\qquad +T\int_0^T\int_{\Omega}|(g_t(x,s+s_i)-g_t(x,s+s_j))w(s)|dxds
\end{align*}
for all $i,j=1,2,\cdots$.
\end{proposition}
{\bf Proof.} Since
\begin{align*}
(g(x,&t+s_i)-g(x,t+s_j))w_t(t) \nonumber
\\
&=\frac{d}{dt}((g(x,t+s_i)-g(x,t+s_j))w(t))-(g_t(x,t+s_i)-g_t(x,t+s_j))w(t),
\end{align*}
we have
\begin{align*}
|\int_0^T\int_s^T&\int_{\Omega}(g(x,\tau+s_i)-g(x,\tau+s_j))w_t(\tau)dxd\tau
ds| \nonumber \\
&\leqslant
\int_0^T\int_{\Omega}|(g(x,T+s_i)-g(x,T+s_j))w(T)|dxds\nonumber
\\
& \qquad+\int_0^T\int_{\Omega}|(g(x,s+s_i)-g(x,s+s_j))w(s)|dxds
\nonumber
\\
& \qquad
+\int_0^T\int_s^T\int_{\Omega}|(g_t(x,\tau+s_i)-g_t(x,\tau+s_j))w(\tau)|dxd\tau
ds.
\end{align*}
Then by $Proposition$ \ref{p5.10} we obtain that for any $T>0$,
\begin{align*}
|\int_0^T\int_s^T&\int_{\Omega}(g(x,\tau+s_i)-g(x,\tau+s_j))w_t(\tau)dxd\tau
ds| \nonumber \\
&\leqslant
2TM\|w(T)\|_{L^2(\Omega)}+2M\int_0^T\|w(s)\|_{L^2(\Omega)}ds
\nonumber
\\
&\qquad +T\int_0^T\int_{\Omega}|(g_t(x,s+s_i)-g_t(x,s+s_j))w(s)|dxds
\nonumber \\
&\leqslant 2TM\|w(T)\|_{L^2(\Omega)}+2MT^{\frac
12}\left(\int_0^T\int_{\Omega}|w(s)|^2dxds\right)^{\frac 12}
\nonumber
\\
&\qquad
+T\int_0^T\int_{\Omega}|(g_t(x,s+s_i)-g_t(x,s+s_j))w(s)|dxds.
\end{align*}
$\hfill \blacksquare$

\begin{proposition}\label{p5.11}
Let $s_i\in \mathbb{R}$ $(i=1,2,\cdots)$, $g \in
L^{\infty}(\mathbb{R};\,L^2(\Omega))\cap
W_{b}^{1,r}(\mathbb{R};\,L^r(\Omega))(r>\frac 65)$,
$\{u_n(t)|~t\geqslant 0,~n=1,2,\cdots\}$ is bounded in
$H_0^1(\Omega)$, and for any $T_1>0$,
$\{u_{n_t}(t)|~~n=1,2,\cdots\}$ is bounded in
$L^{\infty}(0,\,T_1;L^2(\Omega))$. Then for any $T>0$, there exist
subsequences $\{u_{n_k}\}_{k=1}^{\infty}$ of
$\{u_n\}_{n=1}^{\infty}$ and $\{s_{n_k}\}_{k=1}^{\infty}$ of
$\{s_n\}_{n=1}^{\infty}$ such that
\begin{align*}
\lim_{k\to \infty}\lim_{l\to
\infty}\int_0^T\int_s^T\int_{\Omega}(g(x,\tau+s_{n_k})-g(x,\tau+s_{n_l}))(u_{n_k}-u_{n_l})_t(\tau)dxd\tau
ds=0.
\end{align*}
\end{proposition}
{\bf Proof.} Since $\{u_n(t)|~t\geqslant 0,~n=1,2,\cdots\}$ is
bounded in $H_0^1(\Omega)$ and for any $T_1>0$,
$\{u_{n_t}(t)|~~n=1,2,\cdots\}$ is bounded in
$L^{\infty}(0,\,T_1;L^2(\Omega))$, then for any $T>0$, without loss
of generality (at most by passing subsequence), we assume
\begin{equation*}
u_n(T) \to u_0 \quad \text{in}\quad L^2(\Omega)
\end{equation*}
and
\begin{equation*}
u_n \to v \quad \text{in}\quad L^k(0,T;\,L^k(\Omega))~(\text{this
require}~r>\frac 65),
\end{equation*}
where $k<6$.

Then by use of H$\ddot{o}$lder inequality and $Proposition$
\ref{p5.10}, we obtain that
\begin{align*}
&\lim_{n\to \infty}\lim_{m\to
\infty}\int_0^T\int_s^T\int_{\Omega}(g(x,\tau+s_{n})-g(x,\tau+s_{m}))(u_{n}-u_{m})_t(\tau)dxd\tau
ds \nonumber \\
&\leqslant \lim_{n\to \infty}\lim_{m\to
\infty}2MT\left(\int_{\Omega}|u_n(T)-u_m(T)|^2dx\right)^{\frac
12}\nonumber
\\
& \qquad +\lim_{n\to \infty}\lim_{m\to \infty}2MT^{\frac
12}\left(\int_0^T\int_{\Omega}|u_n(s)-u_m(s)|^2dxds\right)^{\frac
12} \nonumber
\\
& \qquad +\lim_{n\to \infty}\lim_{m\to
\infty}\int_0^T\int_{\Omega}|(g_t(x,s+s_n)-g_t(x,s+s_m))(u_n(s)-u_m(s))|dxds
\nonumber \\
&= \lim_{n\to \infty}\lim_{m\to
\infty}T\int_0^T\int_{\Omega}|(g_t(x,s+s_n)-g_t(x,s+s_m))(u_n(s)-u_m(s))|dxds
\nonumber \\
&\leqslant \lim_{n\to \infty}\lim_{m\to
\infty}T\left(\int_0^T\int_{\Omega}|g_t(x,s+s_n)-g_t(x,s+s_m)|^r\right)^{\frac
1r}\left(\int_0^T\int_{\Omega}|u_n(s)-u_m(s)|^k\right)^{\frac 1k}
\nonumber \\
&=0.
\end{align*}
$\hfill \blacksquare$

\subsubsection{A Priori estimates}
The main purpose of this part is to establish
\eqref{5.19}-\eqref{5.21}, which will be used to obtain the
asymptotic compactness.

For any $(u_0^i,\,v_0^i)\in B_0$, let $(u_i(t),\,u_{i_t}(t))$ be the
corresponding solution to $\sigma_i$ with respect to initial data
$(u_0^i,\,v_0^i)$, $i=1,2$, that is, $(u_i(t),\,u_{i_t}(t))$ is the
solution of the following equation
\begin{equation}\label{5.12}
\begin{cases}
u_{tt}+h(u_t)-\Delta u+ f(u(t))=\sigma_i(x,t),\\
u|_{\partial \Omega}=0,\\
(u(0),\,u_t(0))=(u_0^i,\,v_0^i).
\end{cases}
\end{equation}

For convenience, we denote
\begin{equation*}
g_i(t)=\sigma_i(x,\,t),~ h_i(t)=h(u_{i_t}(t)) \quad t\geqslant 0, ~~
i=1,2
\end{equation*}
and
\begin{equation*}
w(t)=u_1(t)-u_2(t).
\end{equation*}

Then $w(t)$ satisfies
\begin{equation}\label{5.13}
\begin{cases}
w_{tt}+h_1(t)-h_2(t)-\Delta w+ f(u_1(t))-f(u_2(t))=g_1(t)-g_2(t),\\
w|_{\partial \Omega}=0,\\
(w(0),\,w_t(0))=(u_0^1,\,v_0^1)-(u_0^2,\,v_0^2).
\end{cases}
\end{equation}

Set
\begin{equation*}
E_w(t)=\frac 12\int_{\Omega}|w(t)|^2+\frac 12 \int_{\Omega}|\nabla
w(t)|^2.
\end{equation*}

{\it Step 1.} Multiplying \eqref{5.13} by $w_t(t)$, and integrating
over $[s,\,T]\times \Omega$, we obtain
\begin{align}\label{5.14}
E_w(T) &+\int_s^T\int_{\Omega}(h_1(\tau)-h_2(\tau))w_t(\tau)dxd\tau
+\int_s^T
\int_{\Omega}(f(u_1(\tau))-f(u_2(\tau)))w_t(\tau)dxd\tau \nonumber \\
&= \int_s^T\int_{\Omega}(g_1(\tau)-g_2(\tau))w_t(\tau)dxd\tau +
E_w(s),
\end{align}
where $0\leqslant s \leqslant T$. Then
\begin{align*}
\int_s^T \int_{\Omega}(h_1(\tau)-h_2(\tau))w_t(\tau)dxd\tau
&\leqslant E_w(s)+
\int_s^T\int_{\Omega}(g_1(\tau)-g_2(\tau))w_t(\tau)dxd\tau \nonumber \\
&\ \quad -
\int_s^T\int_{\Omega}(f(u_1(\tau))-f(u_2(\tau)))w_t(\tau)dxd\tau.
\end{align*}

Combining with $Lemma$ \ref{1}, we get that, for any $\delta >0$,
\begin{align}\label{5.15}
\int_s^T\int_{\Omega}|w_t(\tau)|^2dxd\tau&\leqslant |T-s|\delta
\cdot mes(\Omega) +C_{\delta}E_w(s) + C_{\delta}\int_s^T
\int_{\Omega}(g_1-g_2)w_tdxd\tau
\nonumber \\
&\quad-C_{\delta}\int_s^T
\int_{\Omega}(f(u_1(\tau))-f(u_2(\tau)))w_txd\tau.
\end{align}

{\it Step 2.} Multiplying \eqref{5.13} by $w(t)$ and integrating
over $[0,\,T]\times \Omega$, we get that
\begin{align}\label{5.16}
&\int_0^T\int_{\Omega}|\nabla w(s)|^2dxds + \int_{\Omega}w_t(T)\cdot w(T)dx \nonumber \\
& = \int_0^T\int_{\Omega}|w_t(s)|^2dxds -
\int_0^T\int_{\Omega}(h_1-h_2)wdxds
+ \int_{\Omega}w_t(0)\cdot w(0)dx \nonumber \\
& \qquad - \int_0^T\int_{\Omega}(f(u_1(s))-f(u_2(s)))wdxds+
\int_0^T\int_{\Omega}(g_1-g_2)wdxds.
\end{align}

Therefore, from \eqref{5.15} and \eqref{5.16}, we have
\begin{align*}
2\int_0^TE_w(s)ds & \leqslant 2\delta T mes(\Omega)+2
C_{\delta}E_w(0)+
2C_{\delta}\int_0^T\int_{\Omega}(g_1-g_2)w_tdxds \nonumber \\
& \quad - 2 C_{\delta}\int_0^T\int_{\Omega}(f(u_1(s))-f(u_2(s)))w_t
dxds
- \int_{\Omega}w_t(T)w(T) + \int_{\Omega}w_t(0)w(0) \nonumber \\
& \quad -
\int_0^T\int_{\Omega}(h_1-h_2)w-\int_0^T\int_{\Omega}(f(u_1(s))-f(u_2(s)))w
+\int_0^T\int_{\Omega}(g_1-g_2)w.
\end{align*}

Integrating \eqref{5.14} over $[0,\,T]$ with respect to $s$, we have
that
\begin{align*}
&TE_w(T)+\int_0^T\int_s^T\int_{\Omega}(h_1(\tau)-h_2(\tau))w_t(\tau) dxd\tau ds \nonumber \\
&=-\int_0^T\int_s^T\int_{\Omega}(f(u_1(\tau))-f(u_2(\tau)))w_tdxd\tau ds \nonumber \\
&\quad \qquad+ \int_0^T\int_s^T\int_{\Omega}(g_1-g_2)w_tdxd\tau ds +
\int_0^TE_w(s) ds \nonumber \\
& \leqslant -
\int_0^T\int_s^T\int_{\Omega}(f(u_1(\tau))-f(u_2(\tau)))w_tdxd\tau
ds
+ \int_0^T\int_s^T\int_{\Omega}(g_1-g_2)w_tdxd\tau ds \nonumber \\
& \quad + \delta T mes(\Omega)+C_{\delta}E_w(0) +
C_{\delta}\int_0^T\int_s^T\int_{\Omega}(g_1-g_2)w_tdxds \nonumber
\\
& \quad -C_{\delta}\int_0^T\int_{\Omega}(f(u_1(s))-f(u_2(s)))w_tdxds
-\frac 12 \int_{\Omega}w_t(T)w(T) +\frac 12\int_{\Omega}w_t(0)w(0) \nonumber \\
& \quad -\frac 12 \int_0^T\int_{\Omega}(h_1-h_2)w-\frac12 \int_0^T
\int_{\Omega}(f(u_1(s))-f(u_2(s)))w+ \frac 12
\int_0^T\int_{\Omega}(g_1-g_2)w.
\end{align*}

{\it Step 3.} We need to deal with
$\int_0^T\int_{\Omega}(h_1-h_2)w$. Multiplying \eqref{5.12} by
$u_{i_t}(t)$ we obtain
\begin{align*}
\frac 12\frac{d}{dt}\int_{\Omega}(|u_{i_t}|^2+|\nabla
u_i|^2)+\int_{\Omega}h(u_{i_t})u_{i_t}+\int_{\Omega}f(u_i)u_{i_t}=\int_{\Omega}g_iu_{i_t},
\end{align*}
which, combining with the existence of bounded uniformly absorbing
set, implies that
\begin{align*}
\int_0^T\int_{\Omega}h(u_{i_t})u_{i_t} \leqslant M_T,
\end{align*}
where the constant $M_T$ depends on $T$ (which is different from the
autonomous cases) and the bounds of $B_0$. Then, noticing
\eqref{5.10}, we obtain that
\begin{equation}\label{5.17}
\int_0^T\int_{\Omega}|h(u_{i_t})|^{\frac{p+1}{p}}dxds \leqslant M_T.
\end{equation}

Therefore, using H$\ddot{o}$lder inequality, from \eqref{5.17} we
have
\begin{align*}
|\int_0^T\int_{\Omega}h_iw| & \leqslant
\left(\int_0^T\int_{\Omega}|h(u_{i_t})|^{\frac{p+1}{p}}
\right)^{\frac{p}{p+1}}
\left(\int_0^T\int_{\Omega}|w|^{p+1} \right)^{\frac{1}{p+1}} \nonumber \\
& \leqslant M_T^{\frac{p}{p+1}}\left(\int_0^T\int_{\Omega}|w|^{p+1}
\right)^{\frac{1}{p+1}},
\end{align*}
which implies that
\begin{align}\label{5.18}
|\int_0^T\int_{\Omega}(h_1-h_2)w|  \leqslant
2M_T^{\frac{p}{p+1}}\left(\int_0^T\int_{\Omega}|w|^{p+1}
\right)^{\frac{1}{p+1}}.
\end{align}

\begin{remark}
To some extent, \eqref{5.18} requires that the growth order of $h$
is strictly less  than $5$.
\end{remark}

Set
\begin{align}\label{5.19}
&\phi_{\delta,\,T}((u^1_0,v_0^1),(u^2_0,v_0^2);\sigma_1,\sigma_2) \nonumber \\
& = -
\int_0^T\int_s^T\int_{\Omega}(f(u_1(\tau))-f(u_2(\tau)))w_tdxd\tau
ds \nonumber \\
& \quad + (1+C_{\delta})\int_0^T\int_s^T\int_{\Omega}(g_1-g_2)w_tdxd\tau ds \nonumber \\
& \quad -C_{\delta}\int_0^T\int_{\Omega}(f(u_1(s))-f(u_2(s)))w_tdxds
- \frac 12 \int_0^T\int_{\Omega}(h_1-h_2)wdxds\nonumber \\
& \quad -\frac12 \int_0^T \int_{\Omega}(f(u_1(s))-f(u_2(s)))wdxds+
\frac 12 \int_0^T\int_{\Omega}(g_1-g_2)wdxds,
\end{align}
\begin{equation}\label{5.20}
C_M = \delta T mes(\Omega)+C_{\delta}E_w(0)-\frac 12
\int_{\Omega}w_t(T)w(T)d +\frac 12\int_{\Omega}w_t(0)w(0)dx.
\end{equation}
Then we have
\begin{equation}\label{5.21}
E_w\leqslant\frac{C_M}{T}+\frac1T\phi_{\delta,\,T}((u^1_0,v_0^1),(u^2_0,v_0^2);\sigma_1,\sigma_2).
\end{equation}

\subsubsection{Uniform  asymptotic compactness}

{\bf Proof of Theorem \ref{t5.6}:} \\
Since the family of processes $\{U_\sigma(t,\tau)\}$, $\sigma\in
\Sigma$ has a bounded uniformly absorbing set, for any fixed
$\varepsilon>0$, we can choose first $\delta \leqslant
\frac{\varepsilon}{2mes(\Omega)}$, and then let $T$ so large that
\begin{equation*}
\frac{C_M}{T} \leqslant \varepsilon.
\end{equation*}

Hence, thanks to $Theorem$ \ref{t4.2}, we only need to verify that
$\phi_{\delta,\,T}(\cdot,\,\cdot;\cdot,\,\cdot)\in
\mathfrak{C}(B_0,\,\Sigma)$ for each fixed $T$.

At first, we can observe from the proof procedure of $Theorem$
\ref{t5.5} that for any fixed $T$, we have
\begin{equation}\label{5.22}
\bigcup_{\sigma \in\Sigma}\bigcup_{t \in
[0,\,T]}U_{\sigma}(t,\,0)B_0~~\text{is bounded in}~ E_0,
\end{equation}
and the bound depends on $T$.

At the same time, since $mes_{\mathbb{R}}(\Lambda)=0$, without loss
of generality, we assume $T\notin \Lambda$ (or else, taking
$T_1\notin \Lambda$ satisfies $T_1>T$ and replacing $T$ by $T_1$).

Let $(u_n,u_{n_t})$ be the solutions corresponding to initial data
$(u^n_0,v_0^n)\in B_0$ with respect to symbol $\sigma_n \in \Sigma$,
$n=1,2,\cdots$. Then, from \eqref{5.22}, without loss of generality
(at most by passing subsequence), we assume that
\begin{equation}\label{5.23}
u_n \to u\quad  \star-\text{weakly in} ~ L^{\infty}(0,T;\,
L^6(\Omega)),
\end{equation}
\begin{equation}\label{5.24}
u_n \to u\quad  \text{in} ~L^{p+1}(0,T;\, L^{p+1}(\Omega)),
\end{equation}
\begin{equation}\label{5.25}
u_{n_t} \to u_t\quad  \star-\text{weakly in} ~ L^{\infty}(0,T;\,
L^2(\Omega)),
\end{equation}
\begin{equation}\label{5.26}
u_n \to u\quad  \text{in} ~ L^2(0,T;\, L^2(\Omega))
\end{equation}
and
\begin{equation}\label{5.27}
u_n(0) \to u(0)~~\text{and}~~u_n(T) \to u(T)\quad  \text{in} ~
L^4(\Omega),
\end{equation}
where we use the compact embeddings $H_0^1 \hookrightarrow L^{4}$
and $H_0^1 \hookrightarrow L^{p+1}$\ (since $p<5$).

Now, we will deal with each term corresponding to that in
\eqref{5.19} one by one.

Firstly, from \eqref{5.18} we have
\begin{align*}
|\int_0^T\int_{\Omega}(h(u_{n_t}(s))-&h(u_{m_t}(s)))(u_n(s)-u_m(s))dxds|
\nonumber \\
&\leqslant
2M_T^{\frac{p}{p+1}}\left(\int_0^T\int_{\Omega}|u_n(s)-u_m(s)|^{p+1}
\right)^{\frac{1}{p+1}},
\end{align*}
where $M_T$ depends on $T$ and the norm of $B_0$ in $H_0^1\times
L^2$. Therefore, from \eqref{5.24} we can get
\begin{equation}\label{5.28}
\lim_{n\to \infty}\lim_{m\to
\infty}\int_0^T\int_{\Omega}(h(u_{n_t}(s))-h(u_{m_t}(s)))(u_n(s)-u_m(s))dxds
= 0.
\end{equation}

Secondly, from $Proposition$ \ref{p5.9} and \eqref{5.27}, by the
similar method used in the proof of $Proposition$ \ref{p5.10}, we
can obtain that
\begin{equation}\label{5.29}
\lim_{n\to \infty}\lim_{m\to
\infty}\int_0^T\int_{\Omega}(g_n(x,\,s)-g_m(x,\,s))(u_{n_t}(s)-u_{m_t}(s))dxds=0,
\end{equation}
and from $Proposition$ \ref{p5.11} we can get that
\begin{equation}\label{5.30}
\lim_{n\to \infty}\lim_{m\to
\infty}\int_0^T\int_s^T\int_{\Omega}(g_n(x,\,\tau)-g_m(x,\,\tau))(u_{n_t}(\tau)-u_{m_t}(\tau))dxd\tau
ds=0.
\end{equation}

At the same time, from the growth condition \eqref{1.7} and
\eqref{5.26}, we can get easily that
\begin{equation}\label{5.31}
\lim_{n\to \infty}\lim_{m\to
\infty}\int_0^T\int_{\Omega}(f(u_n(s))-f(u_m(s)))(u_{n}(s)-u_{m}(s))dxds=0.
\end{equation}

Finally, since
\begin{align*}
\int_0^T&\int_{\Omega}
(u_{n_t}(s)-u_{m_t}(s))(f(u_n(s))-f(u_m(s)))dxds \nonumber \\
&=\int_0^T\int_{\Omega} u_{n_t}f(u_n(s))+\int_0^T\int_{\Omega}
u_{m_t}f(u_m(s))
\nonumber \\
&\qquad -\int_0^T\int_{\Omega}
u_{n_t}f(u_m(s))-\int_0^T\int_{\Omega} u_{m_t}f(u_n(s))
\nonumber \\
& =\int_{\Omega} F(u_n(T))-\int_{\Omega} F(u_n(0))+\int_{\Omega}
F(u_m(T))-\int_{\Omega}
F(u_m(0)) \nonumber \\
&\qquad -\int_0^T\int_{\Omega}
u_{n_t}f(u_m(s))-\int_0^T\int_{\Omega} u_{m_t}f(u_n(s)),
\end{align*}
then, by use of \eqref{5.23}, \eqref{5.25}, \eqref{5.27} and
\eqref{1.7}, taking first $m\to \infty$, then $n\to \infty$, we
obtain that
\begin{align}\label{5.32}
\lim_{n\to \infty}&\lim_{m\to \infty}\int_0^T\int_{\Omega}
(u_{n_t}(s)-u_{m_t}(s))(f(u_n(s))-f(u_m(s)))dxds \nonumber \\
&=\int_{\Omega} F(u(T))-\int_{\Omega} F(u(0))+\int_{\Omega}
F(u(T))-\int_{\Omega}
F(u(0)) \nonumber \\
&\qquad -\int_0^T\int_{\Omega}
u_tf(u(s))-\int_0^T\int_{\Omega} u_tf(u(s)) \nonumber \\
&=0.
\end{align}
Similarly, we have
\begin{align*}
\int_s^T&\int_{\Omega}
(u_{n_t}(\tau)-u_{m_t}(\tau))(f(u_n(\tau))-f(u_m(\tau)))dxd\tau \nonumber \\
&=\int_{\Omega} F(u_n(T))-\int_{\Omega} F(u_n(s))+\int_{\Omega}
F(u_m(T))-\int_{\Omega}
F(u_m(s)) \nonumber \\
&\qquad -\int_s^T\int_{\Omega}
u_{n_t}f(u_m(\tau))-\int_s^T\int_{\Omega} u_{m_t}f(u_n(\tau)).
\end{align*}
At the same time, $|\int_s^T\int_{\Omega}
(u_{n_t}(\tau)-u_{m_t}(\tau))(f(u_n(\tau))-f(u_m(\tau)))dxd\tau|$ is
bounded for each fixed $T$, then by Lebesgue dominated convergence
theorem we have
\begin{align}\label{5.33}
&\lim_{n\to \infty}\lim_{m\to \infty}\int_0^T\int_s^T\int_{\Omega}
(u_{n_t}(\tau)-u_{m_t}(\tau))(f(u_n(\tau))-f(u_m(\tau)))dxd\tau ds \nonumber \\
&= \int_0^T\left(\lim_{n\to \infty}\lim_{m\to
\infty}\int_s^T\int_{\Omega}
(u_{n_t}(\tau)-u_{n_t}(\tau))(f(u_n(\tau))-f(u_m(\tau)))dxd\tau
\right)ds \nonumber \\
&=\int_0^T 0ds=0.
\end{align}

Hence, combining \eqref{5.28}-\eqref{5.33}, we get that
$\phi_{\delta,\,T}(\cdot,\,\cdot;\cdot,\,\cdot)\in
\mathfrak{C}(B_0,\,\Sigma)$, and then this completes the proof of
$Theorem$ \ref{t5.6}.           $\hfill \blacksquare$

\begin{remark}\label{r5.13}
If $g_0\in L^{\infty}(\mathbb{R};\,L^2(\Omega))\cap
L_c^2(\mathbb{R};\,L^2(\Omega))$ (e.g., $g_0\in
L^{\infty}(\mathbb{R};\,L^2(\Omega))$ and is a time periodic,
quasi-periodic or almost periodic function in
$L_{loc}^2(\mathbb{R};\,L^2(\Omega))$), then we can obtain
\eqref{5.28} and \eqref{5.29} directly from the definition of
$L_c^2(\mathbb{R};\,L^2(\Omega))$. That is, we do not need all of
the preliminaries from Proposition \ref{p5.8} to Proposition
\ref{p5.11}, and Theorem \ref{t5.6} on uniform asymptotic
compactness still holds.
\end{remark}

\subsection{Existence of uniform attractor}
\begin{theorem}\label{t5.12}
Let $\Omega$ be a bounded domain in $\mathbb{R}^3$ with smooth
boundary, and $h$ and $f$ satisfy \eqref{1.4}-\eqref{1.8}. If
$g_0\in L^{\infty}(\mathbb{R};\,L^2(\Omega)) \cap
W_{b}^{1,\,r}(\mathbb{R};\,L^r(\Omega))$ and $\Sigma$ is defined
by \eqref{5.3}, then the family of processes
$\{U_{\sigma}(t,\tau)\}$, $\sigma\in \Sigma$ corresponding to
\eqref{5.1} or \eqref{1.1} has a compactly uniform (w.r.t.
$\sigma\in \Sigma$) attractor $\mathscr{A}_{\Sigma}$.
\end{theorem}
{\bf Proof.} From $Theorem$ \ref{t5.5} and $Theorem$ \ref{t5.6} we
know that the conditions of $Theorem$ \ref{t3.4} are all satisfied.
$\hfill \blacksquare$

\begin{remark}\label{r5.13a}
If $g_0\in L^{\infty}(\mathbb{R};\,L^2(\Omega))$ and $g_0$ is a time
periodic, quasi-periodic or almost periodic function in
$L_{loc}^2(\mathbb{R};\,L^2(\Omega))$, then the family of processes
$\{U_{\sigma}(t,\tau)\}$, $\sigma\in \Sigma$ corresponding to
\eqref{5.1} or \eqref{1.1} has a compactly uniform (w.r.t.
$\sigma\in \Sigma$) attractor $\mathscr{A}_{\Sigma}$.
\end{remark}

\subsection{Structure of uniform attractor}

In this subsection, we will consider the structure of a uniform
attractor by applying $Theorem$ \ref{t3.8} and $Theorem$
\ref{t3.9}.

For this purpose, we need some continuities for the processes, and
then we need to assume some additional conditions on the external
term $g$, since we need to know whether $\Sigma$ with the $*-$weak
topology of $L^{\infty}(\mathbb{R};L^2(\Omega))\cap
W_b^{1,r}(\mathbb{R};L^2(\Omega))$ is metrizable and when it
becomes a compact metric space.

We assume that $g_0\in W^{1,\infty}(\mathbb{R};L^2(\Omega))$, and
set
\begin{equation}\label{5.34}
\Sigma'_0=\{g_0(x,\,t+h)~|~h\in \mathbb{R}\}
\end{equation}
and
\begin{equation}\label{5.35}
\Sigma'~\text{be the $*-$weakly closure of}~\Sigma'_0 ~\text{in}~
W^{1,\infty}(\mathbb{R};L^2(\Omega)).
\end{equation}

Then, by the classical results (e.g., see \cite{Diestel}), we see
that $\Sigma'$ with the $*-$weak topology of
$W^{1,\infty}(\mathbb{R};L^2(\Omega))$ forms a sequentially compact
and metrizable space. We denote the equivalent metric by
$d_1(\cdot,\,\cdot)$. Thus $(\Sigma',\,d_1)$ is a compact metric
space.

In order to prove the norm-to-weak continuity later, we also need
the following simple property.
\begin{proposition}\label{p5.16}
Let $u_n\rightharpoonup u$ in $L^2(0,T;H_0^1(\Omega))$ and $u_{n_t}
\rightharpoonup u_t$ in $L^2(0,T;L^2(\Omega))$. Then
$u_n(t)\rightharpoonup u(t)$ in $L^2(\Omega)$ for all $t\in
[0,\,T]$.
\end{proposition}
This can be proved by a simple application of $Proposition$ 7.1 and
$Theorem$ 8.1 of Robinson\cite{Ro}, and thus we omit it here.
Moreover, we have the following observation.
\begin{proposition}\label{l6.6.4}
The translation semigroup $\{T(t)\}_{t\geqslant 0}$ is invariant and
continuous in $\Sigma'$ with respect to the $*-$weak topology of
$W^{1,\,\infty}(\mathbb{R};L^2(\Omega))$, equivalently, with respect
to the metric $d_1$.
\end{proposition}\label{p5.17}

\begin{lemma}\label{l5.18}
The family of processes $\{U_\sigma(t,\tau)\}$, $\sigma\in\Sigma'$:
$X\times \Sigma' \mapsto X$ is norm-to-weak continuous.
\end{lemma}

{\bf Proof.} Without loss of generality, we take $\tau =0$.

Multiplying \eqref{5.12} by $u_t$ and integrating over $[0,t]\times
\Omega$,  then from $Theorem$ \ref{t5.5} we can obtain that
\begin{equation}\label{5.36}
\sup_{t\geqslant 0}\|U_{\sigma}(t,0)y_0\|_X\leqslant c(\|y_0\|_X)
\quad \forall~\sigma\in \Sigma', y_0\in X,
\end{equation}
where $c(\cdot)$ is a monotone increasing function on
$\mathbb{R}^+$.

Let $y_{n}\to y_0$ in $X$, $\sigma_n\in \Sigma'$ and $\sigma_n\to
\sigma$ with respect to metric $d_1$. Set
$(u_n(t),\,u_{n_t}(t))=U_{\sigma_n}(t,0)y_{n}$ and
$(u(t),\,u_{t}(t))=U_{\sigma}(t,0)y_{0}$.

Then from \eqref{5.36} we know that $\{U_{\sigma_n}(t,0)y_{n}\}$ is
bounded in $L^{\infty}(0,T;X)$, and due to the growth condition
\eqref{1.7} we also have that $\{f(u_n(t))\}$ is bounded in
$L^{\infty}(0,T;L^2(\Omega))$. At the same time, since $u_n(t)$ is
the solution of \eqref{5.12} we can deduce that $\{u_{n_{tt}}\}$ is
also bounded in $L^{\infty}(0,T; H^{-1}(\Omega))$. Therefore, there
exist subsequence $u_{n_k}$ such that
\begin{equation}\label{5.37}
(u_{n_k}(t),\,u_{n_{k_t}}(t))\to (\tilde{u}(t),\tilde{u}_t(t)) \quad
\text{weakly in}~L^{2}(0,T;X),
\end{equation}
\begin{equation}\label{5.38}
f(u_{n_k})\to \chi \quad *-\text{weakly
in}~L^{\infty}(0,T;L^2(\Omega)),
\end{equation}
\begin{equation}\label{5.39}
u_{n_{k_{tt}}}\to \tilde{u}_{tt}(t) \quad \text{weakly
in}~L^{2}(0,T;\,H^{-1}(\Omega)).
\end{equation}
Similar to that in Lions\cite{Li}, noticing that $u_{n_k}(t)$ is
bounded in $H_0^1(\Omega)$, which implies that $u_{n_k}(t)$ has a
subsequence (here we also denote it by $u_{n_k}(t)$) convergent to
$\tilde{u}(t)$ almost everywhere on $\Omega$, we can get
$\chi=f(\tilde{u}(t))$. Hence, $\tilde{u}(t)$ is a solution
corresponding to initial data $y_0$ with respect to symbol $\sigma$.
Then, by the uniqueness of solution we have that $\tilde{u}=u$.

From \eqref{5.37} and $Proposition$ \ref{p5.16} we have that
$u_{n_k}(t)\rightharpoonup u(t)$ in $L^2(\Omega)$  for all $t\in
[0,T]$. On the other hand, from \eqref{5.36} we know that
$\{u_{n_k}(t)\}$ is bounded in $H_0^1(\Omega)$, therefore,
$u_{n_k}(t)\rightharpoonup u(t)$ in $H_0^1(\Omega)$ for all $t\in
[0,T]$. Hence, for all $t\in [0,T]$, we have
$u_{n}(t)\rightharpoonup u(t)$ in $H_0^1(\Omega)$.

Similarly, from \eqref{5.37}, \eqref{5.39} and the fact that
$\{u_{n_{tt}}\}$ is bounded in $L^{\infty}(0,T; H^{-1}(\Omega))$, we
have that $u_{n_t}(t) \rightharpoonup u_t(t)$ in $L^2(\Omega)$ for
any $t\in [0,T]$. $\hfill \blacksquare$\\

Applying $Theorem$ \ref{t3.9}, from $Theorem$ \ref{t5.5}, $Theorem$
\ref{t5.6}, $Proposition$ \ref{p5.17} and $Lemma$ \ref{l5.18}, we
deduce the following result.

\begin{theorem}\label{t5.19}
Let $\Omega$ be a bounded domain in $\mathbb{R}^3$ with smooth
boundary, and $h$ and $f$ satisfy \eqref{1.4}-\eqref{1.8}. If
$g_0\in W^{1,\,\infty}(\mathbb{R};\,L^2(\Omega))$ and $\Sigma'$ is
defined by \eqref{5.35}, then the family of processes
$\{U_{\sigma}(t,\tau)\}$, $\sigma\in \Sigma'$ corresponding to
\eqref{5.1} has a compactly uniform (w.r.t. $\sigma\in \Sigma'$)
attractor $\mathscr{A}_{\Sigma'}$. Moreover,
\[
\mathscr{A}_{\Sigma'_0}=\mathscr{A}_{\Sigma'}=\bigcup_{\sigma\in
\Sigma'}\mathcal{K}_{\sigma}(0).
\]
\end{theorem}


\end{document}